\documentclass[journal]{IEEEtran}
\usepackage[T1]{fontenc}
\usepackage{algorithm}
\usepackage{amssymb,amsmath}

\usepackage{ifpdf}
\usepackage[pdftex]{graphicx}
\usepackage{array}

\usepackage{url}
\usepackage{color,soul}

\usepackage{subfigure}
\usepackage{float}

\mathchardef\ordinarycolon\mathcode`\:
\mathcode`\:=\string"8000
\begingroup \catcode`\:=\active
  \gdef:{\mathrel{\mathop\ordinarycolon}}
\endgroup

\usepackage{cite}

\ifCLASSINFOpdf
\else
\fi
\usepackage{algorithmic}

\begin{document}

%
\title{A Multi-Grid Iterative Method for Photoacoustic Tomography}
%
%
%

\author{Ashkan~Javaherian~and~Sean~Holman
\thanks{This work was supported in part by a Dean's award from the Faculty of Science and Engineering at the University of Manchester, and the Engineering and Physical Sciences Research Council (EP/M016773/1). \textit{(Corresponding author: Ashkan Javaherian.)}}
\thanks{The authors are with the School
of Mathematics, The University of Manchester, Manchester,
M19 7PL UK (e-mail: ashkan.javaherian@postgrad.manchester.ac.uk).}}

%
%
\markboth{IEEE TRANSACTIONS on MEDICAL IMAGING}%
{Shell \MakeLowercase{\textit{et al.}}}

%

\IEEEpubid{\begin{minipage}{\textwidth}\ \\[6pt] \centering
Copyright (c) 2015 IEEE. Personal use of this material is permitted. However, permission to use this material for any other purposes must be obtained from the IEEE by sending a request to pubs-permissions@ieee.org.
\end{minipage}}


\maketitle

\begin{abstract}
Inspired by the recent advances on minimizing nonsmooth or bound-constrained convex functions on models using varying degrees of fidelity,
we propose a line search multi-grid (MG) method for full-wave iterative image reconstruction in photoacoustic tomography (PAT) in heterogeneous media.
To compute the search direction at each iteration, we decide between the gradient at the target level, or alternatively an approximate error correction at a coarser level, relying on some predefined criteria. To incorporate absorption and dispersion, we derive the analytical adjoint directly from the first-order acoustic wave system. The effectiveness of the proposed method is tested on a total-variation penalized Iterative Shrinkage Thresholding algorithm (ISTA) and its accelerated variant (FISTA), which have been used in many studies of image reconstruction in PAT. The results show the great potential of the proposed method in improving speed of iterative image reconstruction.
\end{abstract}


\section*{Nomenclature}

\addcontentsline{toc}{section}{Nomenclature}
\begin{IEEEdescription}[\IEEEusemathlabelsep\IEEEsetlabelwidth{$V_1,V_2,V_3$}]
\item[MG]Multi-grid.
\item[PAT]Photoacoustic Tomography.
\item[ISTA]Iterative Shrinkage Thresholding Algorithm.
\item[FISTA]Fast Iterative Shrinkage Thresholding Algorithm
\item[FPGA]Field Programmable Gate Array
\item[TR]Time Reversal
\item[GPU]Graphics Processor Unit
\item[AWGN]Additive White Gaussian Noise
\item[TV]Total Variation
\item[PML]Perfectly-Matched Layer
\item[RE]Relative Error
\item[RES]Residual Norm
\item[F]Objective Function (at target level)
\end{IEEEdescription}

\section{Introduction}
%
%
%
%

\IEEEPARstart{P}{hotoacoustic}
Tomography (PAT) is a hybrid imaging technique, which combines the advantage of rich contrast attributed to optical
imaging and high spatial resolution brought up by ultrasound. Typically near-infrared
pulses of light are used to irradiate tissue, which are then absorbed preferentially as a function of the optical absorption of the tissue. The absorbed energy produces local increases in pressure, which move outwards because of the elasticity of soft tissues, and are then sampled temporally by surface detectors \cite{wang-a}.

To estimate the optical absorption distribution of the irradiated tissue from the recorded surface data, one faces two distinct inverse problems, namely acoustic \cite{Treeby} and optical \cite{Tarvainen}. To solve the acoustic inverse problem for media with relatively homogeneous acoustic properties and simple detection surfaces, numerous methods based on filtered back-projection \cite{Kruger,Xu-a,xu-d} or eigenfunction expansion techniques \cite{Xu} have been proposed.

From a practical point of view, real-time 3D reconstruction has been provided for photo-acoustic \cite{Ye} and thermo-acoustic tomography \cite{Ye-b} via 2D reconstruction of slices of the sample, and composing the reconstructed slices into a volume image. Techniques to reduce the effects of out-of-plane acoustic signals were successfully introduced in \cite{Ye}, and in order to achieve the speeds required for real-time imaging reconstruction was done by an onboard FPGA in \cite{Ye-b}. In
these cases, a dense and simple (circular) detection geometry, as well as relatively homogeneous acoustic properties allow
a one-step image reconstruction based on a filtered backprojection algorithm \cite{Zeng}.
\pubidadjcol

TR is a more versatile inversion approach for PAT since it is practical for media with heterogeneous acoustic properties and arbitrary detection geometries \cite{Treeby,Hristova,Hristova-b}. TR and the other inversion approaches mentioned above are inherently based on continuous models, and thus require the detection surface to be very dense and enclose the object \cite{Hristova,Hristova-b}. This is problematic for 3D PAT especially in medical applications.
\IEEEpubidadjcol
The dependance of shape, spectrum and amplitude of propagating acoustic waves on the characteristic properties of tissue media impels enriching the image reconstruction of PAT by simulation of tissue-realistic acoustic propagation \cite{Treeby-e}. To achieve this aim here, the forward problem was solved by a first order acoustic system of three coupled equations which includes two fractional Laplacian operators in order to account separately for absorption and dispersion according to a frequency power law \cite{Treeby-e}.
The main advantage of this acoustic wave propagation model is that it can be efficiently implemented by the k-space pseudospectral method \cite{Cox,Treeby-e} in which the spatial gradient of field parameters is globally computed in frequency domain \cite{Tabei,Treeby-e}.

To mitigate the effects of data incompleteness and noise, iterative methods are often used, e.g., TR-based iterative algorithms \cite{Qian,StefanovUhlmann} or optimization techniques \cite{Wang-c,Huang-c}. Among a great number of optimization approaches, a total variation penalized variant of FISTA \cite{Beck,Beck-b} has been very popular for iterative PAT \cite{Wang-c,Huang-c}. The key element of these optimization approaches is the computation of the gradient of an objective function in terms of the forward model and its adjoint. For heterogeneous media, the adjoint was computed by a ``discretize-then-adjoint'' method in \cite{Huang-c}. Recently, an adjoint was derived for PAT, based on an ``adjoint-then-discretize'' method \cite{Arridge}, from the second order acoustic wave equation which does not include absorption and dispersion. Instead, here the adjoint will be derived using the aforementioned system of three coupled acoustic wave equations \cite{Treeby-e}. In the absence of absorption and dispersion this matches the adjoint in \cite{Arridge}.

In medical PAT, the compartmentalised distribution of chromophores composing tissues induces step-like pressure discontinuities within absorbing regions, which make the generated waves highly broadband \cite{Treeby}. To cover such a broad range of frequencies in the reconstruction, very dense grids are needed which make iterative PAT computationally burdensome. To mitigate this problem, numerous methods have been proposed to accelerate wave propagation models\cite{Dean-Ben,Lutzweiler-b,Arridge-b}. By recent advances on data casting and parallelization using GPUs \cite{Dean-Ben-b,Wang-b,Treeby-b} or FPGA-based hardware implementation of the reconstruction algorithms \cite{Ye-b}, wave propagation models were accelerated notably. In the present work, we take a different approach, and look at a method to improve the performance of the underlying algorithm.

In general, whenever a finite-dimensional optimization problem arises from an infinite-dimensional continuous problem, it is possible to control the fidelity with which the optimization model captures the underlying continuous problem \cite{Parpas}. In the case that the arising discretized model is very large-scale, multi-grid (MG) schemes, which exploit a hierarchy of discretized models (levels) of varying size, are very popular.

A MG scheme for unconstrained smooth optimization problems was first proposed by Nash \cite{Nash}, in which the information at the coarse level is utilized to compute the search direction at the target level. This method was recently extended to composite convex functions involving a smooth term plus a non-smooth $\ell_1$ term \cite{Parpas,Hovhannisyan}, relying on the recent theoretical advances on minimization of smoothable functions \cite{Beck-c}. Additionally, Nash's method was recently extended to smooth bound-constrained optimization \cite{Kocvara}.

To mitigate the burdensome computational requirements of iterative PAT, we propose a line search multi-grid method for full-wave iterative image reconstruction in PAT so that at some iterations, a recursive search direction is computed by minimizing the objective function at some coarser levels. Here the proposed MG method is applied to ISTA and FISTA on two levels, but it can be easily extended to other first order methods such as Primal-Dual algorithms, or to more than two levels.

\section{Background}

\subsection{Forward Problem}      \label{aa}
\subsubsection{Lossless media}
Acoustic propagation in lossless heterogeneous media can be described relying on three coupled equations, i.e, equation of motion, equation of continuity, and adiabatic equation of state, respectively in the form \cite{Treeby}
\begin{align}\label{mom}
\frac{\partial u}{\partial t}(r,t) & = - \frac{1}{\rho_0(r)} \nabla p(r,t),\\
\label{mass}
\frac{\partial \rho}{\partial t}(r,t) & = - \rho_0(r)\ \nabla \cdot u(r,t),\\
\label{estate}
p(r,t) & = c_0(r)^2  \rho (r,t)
\end{align}
with initial conditions
\begin{equation} \label{initial}
p(r,0) = p_0(r), \quad u(r,0) = 0.
\end{equation}
Here, $p(r,t)$ denotes the acoustic pressure at position $r\in\mathbb{R}^d$ ($d=2 \hspace{0.1cm}\text{or}\hspace{0.1cm}3$), and time $t\in\mathbb{R}^{+} $. Additionally, $u(r,t)$ denotes the vector-valued acoustic particle velocity, $c_0(r)$ denotes the varying sound speed, and $\rho(r,t)$ and $\rho_0(r)$ represent the acoustic and ambient densities, respectively.

\subsubsection{Lossy media}
To simulate wave propagation in lossy media, a trade-off is typically needed between agreement with experimental observations \cite{Szabo}, meeting causality conditions \cite{Nachman,Kowar,Kowar-b}, and efficiency of numerical computations \cite{Chen,Treeby-e}. Over frequencies of generated ultrasound waves in PAT, the absorption in tissue obeys a frequency power law in the form
\begin{equation}
\alpha=\alpha_0 w^y,
\end{equation}
where $\alpha_0$ is the absorption coefficient in $\text{dB} \hspace{0.1cm} \text{MHz}^{-y} \hspace{0.1cm} \text{cm}^{-1}$, $w$ is the angular frequency in $\text{MHz}$, and $y$ is the power law exponent \cite{Treeby}.

Describing the wave attenuation as effects of viscosity and thermal conduction leads to the so-called thermo-viscous attenuation model, which yields a frequency-squared attenuation ($y=2$). This model does not match the observed frequency dependence of attenuation in tissues \cite{Szabo}, and also violates the causality condition \cite{Kowar}. The attenuation model was later described by a superposition of relaxation mechanisms \cite{Nachman}. This model meets the causality condition \cite{Kowar}, and can be implemented efficiently by the k-space method  \cite{Tabei}. However, for simulating broadband acoustic propagation in PAT, estimation of the distribution of relaxation parameters for each relaxation process is troublesome. To account for the power law dependence on frequency evident in biological tissue ($1<y<1.5$), a lossy wave equation based on temporal convolution, the so-called Szabo's model \cite{Szabo}, was proposed, which was later rewritten as a time-domain fractional derivative operator, e.g., \cite{Liebler,Chen-b}. Szabos's model has been shown to be noncausal for $y>1$ \cite{Kowar}. Furthermore, the time-domain convolution or fractional derivative operators inherently require storing the complete pressure field at previous times. A memory-efficient power law absorption model based on fractional Laplacian operators was proposed in \cite{Chen}, and was then modified to incorporate the dispersive sound speed \cite{Treeby-e}. This model can be easily incorporated into the k-space pseudospectral method without storing the computed pressure field at previous time steps \cite{Treeby-e}, as opposed to classical absorption and dispersion models that involve time-domain fractional operators. This model is very popular in PAT \cite{Treeby,Treeby-b,Treeby-e}, and was thus used in the present study.

Applying this model, the absorption and dispersion effects are incorporated into the wave propagation by adding two fractional Laplacian operators to the equation of state in the form \cite{Treeby-e}
\begin{equation}\label{state}
\begin{split}
p(r,t) = c_0(r)^2 \Big \{& 1 - \tau(r) \frac{\partial}{\partial t} (-\nabla^2)^{\frac{y}{2}-1}\\
&- \eta(r) (-\nabla^2)^{\frac{y-1}{2}} \Big \} \rho (r,t).
\end{split}
\end{equation}
The absorption and dispersion proportionality coefficients, $\tau(r)$ and $\eta(r)$ respectively, are calculated by
\begin{align}
&\tau(r)=-2\alpha_0 c_0(r)^{y-1}, \quad \eta(r)=2\alpha_0 c_0(r)^y \tan(\pi y/2).     \label{Twelve}
\end{align}
A toolbox for modeling of acoustic wavefield propagation based on the k-space pseudo-spectral method is freely available \cite{Treeby-b}, and was used in this study.

Recently it has been shown that this attenuation model encounters some noncausality problems, since the corresponding equation of state \eqref{state} is nonlocal in space at each time instant, and in addition the Green's function of the resulting wave equations does not have a finite wave front speed \cite{Kowar-b}. To maintain these causality conditions, the state equation was spatially localized by enforcing a local time shift to the attenuation model. This also leads to a Green's function with a finite wave front speed, which is known to be a strong causality condition for systems of wave equations \cite{Kowar,Kowar-b}. However, similar to classical lossy wave equations based on temporal fractional derivatives \cite{Liebler,Chen-b}, the numerical implementation of this model is very expensive for 3D PAT.

\subsection{Inverse Problem}    \label{bb}
The acoustic inverse problem is to estimate the initial pressure $p_0(r)$ inside a bounded region $\Omega \subset \mathbb{R}^d$ from the measurements $p_m(r_s,t)$ taken at positions $r_s$ within an open set $\Gamma \subset \mathbb{R}^{d-1}$ from time $t=0$ to $T$.

\subsubsection{Time Reversal (TR)}
Employing the time reversal method, $p_m(r_s,t)$ is enforced as a Dirichlet boundary condition in a time-reversed order, yielding
\begin{align}
p_\text{tr}(r_s,t)=p_m(r_s,T-t).
\end{align}
Here, the time $T$ is assumed to be sufficiently large so that all waves leave the medium, yielding $p_\text{tr}(r,0)=0$ \cite{Treeby}. This is, however, not exactly held for even dimensions or heterogeneous media\cite{Hristova,Hristova-b}. It was shown that to account for the absorption and physical dispersion in TR,
the absorption term in \eqref{state} must be reversed in sign, while the physical
dispersion term remains unchanged \cite{Treeby}.

\subsubsection{Variational methods}
The accuracy of inversion approaches, including TR, is limited for sub-sampled or noisy data.
In these cases, variational image reconstruction methods provide an effective alternative \cite{Arridge,Huang-c}. Let $p_0$ denote the sought after
initial pressure distribution, and $\hat{p} \in \mathbb{R}^{M} (M \in \mathbb{N})$ and $\varepsilon$ denote the time series of measured data at sensors and the corresponding noise, respectively. Additionally, let $\mathcal{H}$ represent the forward model discussed in \ref{aa}. We then have
\begin{align}
\hat{p}=\mathcal{H}p_0+\varepsilon.
\end{align}
Inferring $p_0$ from $\hat{p}$ amounts to solving a regularized least-square optimization problem in the form
\begin{align}
\providecommand{\norm}[1]{\left\lVert#1\right\rVert}
p = \operatornamewithlimits{argmin}\limits_{p_0\geqslant{0}}\frac{1}{2} \norm{\mathcal{H} p_0-\hat{p}}^2+ \lambda \mathcal{J}(p_0).  \label{One}
\end{align}
\noindent Here, $\lambda>0$ is a regularization parameter, and $\mathcal{J}(p_0)$ is a regularization functional that can be used to impose a-priori information
about the true solution. Here, regularization functional $\mathcal{J}$ is taken to be total variation (TV) since it is very popular in PAT because of accounting for feature edges in the reconstruction \cite{Arridge-b,Huang-c}.

Solving \eqref{One} requires the computation of gradient of the objective function as a function of the forward operator $\mathcal{H}$ and its adjoint $\mathcal{H}^*$.
Specifically in PAT, the main computational cost of the minimization problem is the implementation of $\mathcal{H}$ and $\mathcal{H}^*$, and the cost of other steps is negligible in comparison. To derive the adjoint for the forward model described in \ref{aa}, a ``discretize-then-adjoint'' method was proposed, where the computational steps of the discretized forward problem are explicitly reversed \cite{Huang-c}. The adjoint obtained by this strategy may not correspond exactly with a discretization of the adjoint in the continuous domain. Very recently, a general analytic form of the adjoint in PAT was derived \cite{Arridge}, where the time-reversed pressure $p_m(r_s,T-t)$ is added as a time-dependent mass source term $s(r,t)$ to \eqref{mass}. In comparison, in the TR approach this is enforced as an explicit Dirichlet boundary condition \cite{Treeby}. However, the method in \cite{Arridge} gives the adjoint using the second order acoustic equation, which does not include absorption and dispersion. In order to include absorption and dispersion, we will derive the adjoint using \eqref{mom}, \eqref{mass} and \eqref{state}.

\section{Adjoint of the Three-Coupled First Order Wave Propagation Equation}   \label{A1}
\maketitle
Similar to \cite{Arridge}, the continuous forward operator is the map
\begin{equation}
\begin{split}
&\mathcal{H}: C_0^\infty(\Omega) \rightarrow \mathbb {R}^M,\\
&\mathcal{H} [p_0] (r,t)= \mathcal{M}  w(r,t) p(r,t),
\end{split}
\end{equation}
where $w(r,t) \in C_0^\infty (\Gamma \times [0,T])$ restricts the pressure $p(r,t)$ to the spatio-temporal field accessible to the sensors, and $\mathcal{M}$ maps the accessible part of the pressure field into the measured data at sensors $\hat{p} \in \mathbb{R}^M$. Like \cite{Arridge}, we will assume that $\mathcal{M}^*$ is given. Now, we have $\mathcal{H}^*=\mathcal{P}^*\mathcal{M}^*$, where $\mathcal{P}^*:C_0^\infty(\Gamma \times [0,T])\rightarrow C_0^\infty (\Omega)$ is the adjoint of
\begin{equation}
\begin{split}
&\mathcal{P}: C_0^\infty(\Omega) \rightarrow    C_0^\infty(\Gamma \times [0,T])  \\
&\mathcal{P} [p_0] (r,t)=   w(r,t) p(r,t).
\end{split}
\end{equation}
Let us first define the time-reversed adjoint fields, $p^*$, $u^*$, $\rho^*$ by
\begin{align}
\label{moma}
\frac{\partial u^*}{\partial t}(r,t)  & = - \frac{1}{\rho_0(r)} \nabla p^*(r,t),\\
\label{massa}
\frac{\partial \rho^*}{\partial t}(r,t) & =  \rho_0(r) \Big ( - \nabla \cdot u^*(r,t) \\
&\quad\quad\quad + w(r,T-t) h(r,T-t) \Big ), \nonumber\\
\label{statea}
p^*(r,t) & =  \rho_0(r) \Big \{ 1 - \frac{\partial}{\partial t} (-\nabla^2)^{\frac{y}{2}-1} \tau(r)\\
&\quad\quad\quad - (-\nabla^2)^{\frac{y-1}{2}}  \eta(r)  \Big \} \frac{c_0(r)^2}{\rho_0(r)} \rho^* (r,t) \nonumber
\end{align}
with initial conditions
\begin{equation} \label{initial2}
p^*(r,0) = 0, \quad u^*(r,0) = 0.
\end{equation}
By definition of the adjoint, for any $h(r,t) \in C_0^\infty (\Gamma \times [0,T])$
\begin{equation}\label{adj2}
\int_0^T \int_{\mathbb{R}^d}  \mathcal{P} [p_0](r,t) \ h(r,t) \ \mathrm{d} r\ \mathrm{d} t = \int_{\mathbb{R}^d} p_0(r)\  \mathcal{P}^* [h](r) \ \mathrm{d} r.
\end{equation}
The claim is that
\begin{equation} \label{adjj}
\mathcal{P}^* [h](r) = \frac{\rho^*(r,T)}{\rho_0(r)}.
\end{equation}
To prove this we start with \eqref{mass}, which yields
\begin{equation}
\begin{split}
0  = & \int_0^T  \int_{\mathbb{R}^d}  \left ( \frac{\partial \rho}{\partial t} + \rho_0 \nabla \cdot u \right ) \frac{p^*(r,T-t)}{\rho_0}  \mathrm{d} r \ \mathrm{d} t. \\
  = & \int_0^T  \int_{\mathbb{R}^d}  \left ( \frac{\partial \rho}{\partial t} + \rho_0 \nabla \cdot u \right ) \Big \{ 1 - \frac{\partial}{\partial t} (-\nabla^2)^{\frac{y}{2}-1} \tau\\
  &\hskip1cm - (-\nabla^2)^{\frac{y-1}{2}}  \eta  \Big \} c_0^2 \frac{\rho^* (r,T-t)}{\rho_0} \ \mathrm{d} r \ \mathrm{d} t.
 \end{split}
\end{equation}
where we also used \eqref{statea} for the second equality, and suppressed the dependence on $(r,t)$ in some places for brevity. Now, the Laplacian $-\nabla^2$ is self-adjoint, and so the fractional powers of the Laplacian are as well and we can move them from acting on $c_0^2 \rho^*/\rho_0$ to the first term as in integration-by-parts. This gives
\begin{equation}
\begin{split}
0 = & \int_0^T \int_{\mathbb{R}^d} \left [ \frac{\partial}{\partial t}c_0^2 \Big \{ 1 - \eta (-\nabla^2)^{\frac{y-1}{2}} \Big \} \rho\right ]  \frac{\rho^*(r,T-t)}{\rho_0} \\
    &\hskip1cm - \left [c_0^2 \tau \frac{\partial}{\partial t} (-\nabla^2)^{\frac{y}{2}-1} \rho \right ] \frac{\frac{\partial \rho^*}{\partial t}(r,T-t)}{\rho_0} \\
    &\hskip2cm + \nabla \cdot u \ p^*(r,T-t) \ \mathrm{d}r\ \mathrm{d} t.
\end{split}
\end{equation}
Now we apply integration-by-parts to the terms on the first and third lines as well as \eqref{state}, and the initial conditions in \eqref{initial2}, which together yield
\begin{equation}
\begin{split}
0 = & \int_0^T  \int_{\mathbb{R}^d} \frac{p}{\rho_0}  \frac{\partial \rho^*}{\partial t}(r,T-t) - u \cdot \nabla p^*(r,T-t) \ \mathrm{d}r\ \mathrm{d} t\\
    & - \int_{\mathbb{R}^d} \left [c_0^2 \Big \{ 1 - \eta (-\nabla^2)^{\frac{y-1}{2}} \Big \} \rho(r,0)\right ] \frac{\rho^*(r,T)}{\rho_0} \ \mathrm{d} r.
\end{split}
\end{equation}
Considering that $\frac{\partial \rho}{\partial t}(r,0) = -\rho_0(r) \nabla \cdot u(r,0) = 0$ by \eqref{mass}, and in light of \eqref{state}, we have
\begin{equation}
c_0^2 \Big \{ 1 - \eta (-\nabla^2)^{\frac{y-1}{2}} \Big \} \rho(r,0) = p_0(r).
\end{equation}
Putting this into the previous formula finally gives
\begin{equation}\label{gh1}
\begin{split}
&\int_{\mathbb{R}^d} p_0(r) \frac{\rho^*(r,T)}{\rho_0(r)} \ \mathrm{d} r =\\
&\int_0^T  \int_{\mathbb{R}^d} \frac{p}{\rho_0}  \frac{\partial \rho^*}{\partial t}(r,T-t) - u \cdot \nabla p^*(r,T-t) \ \mathrm{d}r\ \mathrm{d} t.
\end{split}
\end{equation}
Putting \eqref{gh1} aside for a moment, we next use \eqref{mom} which gives
\begin{equation}
\begin{split}
0 & = \int_0^T \int_{\mathbb{R}^d}  \left (\rho_0 \frac{\partial u}{\partial t} + \nabla p\right ) \cdot u^*(r,T-t) \ \mathrm{d} r \ \mathrm{d} t.
\end{split}
\end{equation}
Applying integration-by-parts to both terms, and enforcing the initial conditions \eqref{initial} and \eqref{initial2} gives
\begin{equation}\label{gh2}
0 = \int_0^T \int_{\mathbb{R}^d}  \rho_0 \ u \cdot \frac{\partial u^*}{\partial t}(r,T-t) -  p\ \nabla \cdot u^*(r,T-t) \ \mathrm{d} r \ \mathrm{d} t.
\end{equation}
Subtracting \eqref{gh2} from \eqref{gh1} and gathering on the right-hand-side the terms involving $p$ and $u$ respectively in lines two and three of the next formula, we have
\begin{equation}
\begin{split}
&  \int_{\mathbb{R}^d} p_0(r) \frac{\rho^*(r,T)}{\rho_0(r)} \ \mathrm{d} r  = \\
&\int_0^T \int_{\mathbb{R}^d}  \left (\frac{1}{\rho_0} \frac{\partial \rho^*}{\partial t}(r,T-t) + \nabla \cdot u^*(r,T-t) \right ) p \ \mathrm{d} r\ \mathrm{d} t\\
&- \int_0^T \int_{\mathbb{R}^d}  \left (\rho_0 \frac{\partial u^*}{\partial t}(r,T-t) + \nabla p^*(r,T-t) \right ) \cdot u  \ \mathrm{d} r\ \mathrm{d} t.
\end{split}
\end{equation}
Now, by \eqref{massa} the integrand in the second line above is equal to $p(r,t)w(r,t)h(r,t)$, and by \eqref{moma} the third line is equal to zero. Therefore
\begin{equation}
\begin{split}
\int_{\mathbb{R}^d} p_0(r) \frac{\rho^*(r,T)}{\rho_0(r)} \ \mathrm{d} r & =  \int_0^T \int_{\mathbb{R}^d} w(r,t)  p(r,t) h(r,t)  \ \mathrm{d} r \ \mathrm{d} t \\
&= \int_0^T \int_{\mathbb{R}^d}  \mathcal{P} [p_0](r,t) \ h(r,t) \ \mathrm{d} r\ \mathrm{d} t.
\end{split}
\end{equation}
Finally, using \eqref{adj2}, we see now that $\mathcal{P}^* [h](r) = \frac{\rho^*(r,T)}{\rho_0(r)}$, and thus the claim \eqref{adjj} about the adjoint is proven. Taking $\tau(r)=0$ and $\eta(r)=0$ makes this adjoint the same as that proposed in \cite{Arridge} for lossless media.

\section{First-order Optimization Methods for PAT}    \label{Thirty-six}
The numerical implementation of the derived forward and adjoint operators requires discretization of the models.
Accordingly, the discretized variant of the sought after initial pressure $p_0$ is denoted by $x \in \mathbb{R}^N (N \in \mathbb{N})$ with $N$ the number of grid points, and
the discretized forward model linking $x$ to data $\hat{p}$ is denoted by $H \in \mathbb{R}^ {M \times N}$.
Problem \eqref{One} is in a class of non-smooth constrained convex minimization problems of the form
\begin{align}
\operatornamewithlimits{argmin}\limits_{x}\left\{F(x):=f(x)+g(x)\right\}.   \label{two}
\end{align}
Here,
$\providecommand{\norm}[1]{\left\lVert#1\right\rVert}
f(x)= \frac{1}{2}\norm{Hx-\hat{p}}^2$ is a continuously differentiable function with Lipschitz continuous gradient having smallest Lipschitz constant $\providecommand{\norm}[1]{\left\lVert#1\right\rVert}
L_f=\sigma_{\text{max}}(H^*H)$, where $\sigma_{\text{max}}(.)$ stands for the largest singular value. The gradient of $f$ is computed by
\begin{align}
\nabla{f}(x)=H^*\left(Hx-\hat{p}\right).    \label{Six}
\end{align}
Additionally, we take
$g(x)= \lambda \mathcal{J}(x)+\delta_C\left({x}\right)$ where $\delta_C$ is the indicator function for the set of constraints $C = \{x \geqslant 0\}$.

Applying the so-called forward-backward splitting method to a fixed point iterative scheme arising from the optimality conditions of problem \eqref{two} gives the two-step Iterative Shrinkage Thresholding Algorithm (ISTA) shown in Algorithm \ref{ISTA}.
\begin{algorithm}
\caption {ISTA}
\begin{algorithmic}[1] \label{ISTA}
\STATE {\textbf{Iteration 0:} $x_0$}
\STATE {\textbf{Iteration $\mathbf{k \geqslant 1}$:}}
\STATE  $ z_k = x_{k-1}-\alpha_k\ \nabla{f}(x_{k-1})$   \label{smooth}
\STATE  $x_k = \text{prox}_{\alpha_k}(g)(z_k)$ \label{proximal}
\STATE   {\textbf{Output:}   $x_*$}.
\end{algorithmic}
\end{algorithm}

\noindent Here, line \ref{smooth} is a forward gradient descent step \cite{Goldstein}, and is in a class of line search techniques which utilize a steepest descent search direction $- \nabla f(y)$ and step size $\alpha_k$ \cite{Beck}. Applying ISTA, the convergence of the iterates $x_k$ to a minimizer $x_*$ of problem \eqref{two} is proven if $\alpha_k \in \left(0, 2/L_f\right)$ \cite{Beck-b}. To determine $L_f$, the largest singular value of $H^*H$ is computed iteratively by the power method \cite{Arridge}. Since $L_f$ is independent to the unknown $x$, it can be stored and used for all experiments done in a fixed setting \cite{Arridge-b}. Otherwise, $\alpha_k$ can be computed adaptively by backtracking line search techniques, although this is inefficient for large-scale PAT problems. Additionally,
\begin{align}
\providecommand{\norm}[1]{\left\lVert#1\right\rVert}
\begin{split}
     \text{prox}_{\alpha_k}(g)(z_k):=
     \operatornamewithlimits{argmin}\limits_{x}\left\{ g(x)+\frac{1}{2\alpha_k}\norm{x-z_k}^2 \right\}.
\label{proxmal}
\end{split}
\end{align}
is a backward gradient step, and is called the proximal map \cite{Goldstein}. Similar to \cite{Huang-c}, here the proximal map associated with the TV functional was computed based on a dual approach given in \cite{Beck-b}.

The computational cost of performing the forward and adjoint solvers necessary to compute $\nabla f$ on a grid of size $N_x \cdot N_y \cdot N_z$ in $N_t$ time steps is $\mathcal{O}(N_t \cdot N_x \cdot N_y \cdot N_z  \log (N_x \cdot N_y \cdot N_z))$, whereas the computational cost of solving the proximal map by the dual approach given in \cite{Beck-b} is $\mathcal{O}(N_x \cdot N_y \cdot N_z)$. As a result, the major cost of each iteration is the forward gradient step, while the cost of the proximal map in \eqref{proxmal} is almost negligible in comparison.

An acceleration to ISTA is given by FISTA \cite{Beck}, which provides a global convergence rate of $\mathcal O(1/\sqrt {\epsilon})$, compared to $\mathcal O(1/\epsilon)$ for ISTA, where $\epsilon$ denotes the desired accuracy. FISTA is outlined in Algorithm \ref{FISTA} \cite{Beck}.

\begin{algorithm}
\caption {FISTA}
\begin{algorithmic}[1] \label{FISTA}
\STATE {\textbf{Iteration 0:} $y_1=x_0, t_k=1 \quad (\theta_k=0)$}
\STATE {\textbf{Iteration $\mathbf{k \geqslant 1}$:}}
\STATE  $x_k =P_L(y_k)$   \label{pro}
\STATE  $t_{k+1} = \frac{1+\sqrt{1+4{t_k}^2}}{2}$
\STATE  $\theta_k = \frac{t_{k}-1}{t_{k+1}}$ \label{thetak}
\STATE  $y_{k+1}  = x_k+ \theta_k \left(x_k-x_{k-1}\right)$
\STATE   {\textbf{Output:}   $x_*$}
\end{algorithmic}
\end{algorithm}
\noindent Here, operator $P_L(\cdot)$ represents lines \ref{smooth} and \ref{proximal} in Algorithm \ref{ISTA}.
Note that replacing line \ref{thetak} by $\theta_k=0$ in Algorithm \ref{FISTA} gives $y_{k+1}=x_k$, and reduces the algorithm to ISTA. In the next section, the multi-grid algorithm will be described for a general algorithm like FISTA but with line \ref{thetak} possibly replaced. Thus an extension to ISTA or other first-order optimization methods is straightforward. The convergence of sequence $x_k$ provided by FISTA is proven when $\alpha_k \in \left(0, 1/L_f\right)$ \cite{Beck}. For applications of ISTA in iterative PAT, see \cite{Arridge}, and for FISTA, see \cite{Huang-c,Wang-c}.

\section{Line search Multi-grid Optimization Method}
To improve the speed of Algorithm \ref{FISTA}, a multi-grid (MG) line search strategy is adopted based on Nash's well-known method \cite{Nash} so that at each iteration the algorithm decides between two possibilities: a direct search direction computed at the target level, or alternatively a recursive search direction generated from some steps taken at coarser levels.

Considering the computational cost of the forward and adjoint operators given in Section \ref{Thirty-six}, the cost of performing them on a coarse grid with a size $(N_x/2) \cdot (N_y/2) \cdot (N_z/2)$ in $N_t/2$ time steps is less than $1/16$ the cost on the fine grid, as the time step is changed proportionally to the spatial distance of grid points. Accordingly, for 2D PAT, the computational cost of coarse forward and adjoint models is less than $1/8$ the cost on the fine model.

\subsection{First-order Coherence of Levels for Smooth Unconstrained Optimization: An Extension to FISTA}   \label{coherence}
We denote the level that supports the fine resolution, referred to here as the ``target level'', by subscript $h$ and the next coarse level by $h-1$.
The transfer of information from level $h$ to $h-1$ is done by restriction operator $I_h^{h-1}$. Conversely, prolongation operator $I_{h-1}^h$ is used to transfer information from level $h-1$ to $h$.

To guarantee convergence on multiple levels, the first order optimality conditions of the levels must match.
To attain this, Nash \cite{Nash} suggests adding a linear term to the objective function at the next coarse level. We extend Nash's method to FISTA so that to compute a recursive search direction, starting from iteration $y_{h,k}$ at the target level, we use as the objective function at the next coarse level $h-1$
\begin{align}
\phi_{h-1}\left(x_{h-1}\right)=F_{h-1}\left(x_{h-1}\right)+\langle v_{h-1},x_{h-1} \rangle,    \label{fourteen}
\end{align}
where $v_{h-1}$ stands for
\begin{align}
v_{h-1}= I_h^{h-1} \nabla F_h \left(y_{h,k}\right) - \nabla F_{h-1} \left(x_{h-1,0}\right),     \label{Fourteen}
\end{align}
with ${x_{h-1,0}}=I_h^{h-1} y_{h,k}$ the initial point at the next coarse level $h-1$. Note that $y_{h-1,1}=x_{h-1,0}$ according to the initialization in Algorithm \ref{FISTA}. In this way, the gradient of the objective functions at the point of transfer between the two levels matches so that
\begin{align}
\nabla \phi_{h-1} \left(x_{h-1,0}\right)=I_h^{h-1} \nabla F_h \left(y_{h,k}\right).
\end{align}
This property is called ``first-order coherence'' \cite{Nash,Wen,Parpas,Hovhannisyan}.

\subsection{Extension to Non-smooth Unconstrained Optimization}
The approach given above is not applicable to non-smooth objective functions since the computation of $\nabla F$ is not possible.
From a theoretical point of view, an approach for minimizing non-smooth functions via treating the problem as a sequence of smooth problems has been considered. A global convergence rate of $\mathcal O(1/\epsilon)$ was first established for functions with so-called ``explicit max-structure" \cite{Nesterov}, and was then extended to the so-called ``smoothable" functions \cite{Beck-c}.

Recently, relying on the mentioned works, Nash's multi-grid method was extended to unconstrained composite functions involving a smooth term plus a nonsmooth $\ell_1$ term. The convergence rate of ISTA on a multi-grid setting was established by \cite{Parpas}. Recently, an MG method with an optimal rate of convergence ($\mathcal O(1/\sqrt{\epsilon})$) was proposed \cite{Hovhannisyan}. This has been inspired by a modified variant of Nesterov's acceleration technique, where the problem is treated as a linear combination of primal gradient and mirror descent steps \cite{Allen-Zhu}. The global convergence rate established by this MG method is optimal, but the bound on the
worst case convergence rate is greater than that of the standard ``gradient and mirror descent" algorithm in \cite{Allen-Zhu}. Note that in practice the sequence $x_k$ provided by the ``gradient and mirror descent'' algorithm matches that of FISTA on a fixed grid \cite{Allen-Zhu,Hovhannisyan}. However, as opposed to the MG variant of FISTA proposed here, we observed that the MG algorithm proposed in \cite{Hovhannisyan} is not efficient in PAT.

In order to use MG with FISTA, at $y_{h,k}$, we smooth the TV penalty function in the form
\begin{equation}
\providecommand{\abs}[1]{\lvert#1\rvert}
\mathcal J_\rho(y_{h,k})=\sum_{{n_1,n_2,n_3}} \sqrt {{\abs{(\nabla{y_{h,k}})_{n_1,n_2,n_3}}}^2+\rho^2}- \rho,   \label{Thirty-two}
\end{equation}
where $\rho$ is the smoothing parameter. The gradient of $F$
is now computed as \cite{Vogel}
\begin{equation}
\providecommand{\abs}[1]{\lvert#1\rvert}
\nabla F_{\rho}(y_{h,k})=H^*\left(Hy_{h,k}-\hat{p}\right)\\
- \lambda \hspace{0.1cm} \nabla \cdot \left(\frac{\nabla y_{h,k}}{\sqrt {\abs{\nabla y_{h,k}}^2+\rho^2}}\right).   \label{Twentyone}
\end{equation}
The implementation of Nash's method via computing $v_{h-1}$ and minimizing $\phi_{h-1}$ by Eqs. (\ref{Fourteen}) and (\ref{fourteen}) is now straightforward in the unconstrained case.

\subsection{Extension to Constrained Convex Optimization}
The coherence formula does not account for the bound constraint that is enforced in the PAT problem.
In general, very few studies exist in the literature to extend Nash's method to bound-constrained optimization, e.g., \cite{Gratton}.
In \cite{Graser}, a method to deal with bound constraints for MG optimization problems was proposed based on truncation of the set of indices at which the constraints are active. In this method, the active nodes on the fine level are fixed for the next coarse-grid correction. This truncation scheme is very conservative, and thus makes the MG algorithm inefficient.

Recently an MG approach for smooth constrained optimization problems has been developed via restriction of bound constraints, rather than the truncation of active set \cite{Kocvara}. The restriction of constraints is done so that a feasible point remains feasible after the coarse correction step. In our specific case, this MG approach is applied to two levels, and a nonnegativity constraint is enforced globally to all nodes at the fine level. Let $\mathcal I_{h,i}$ denote the union of indices at level $h$ that locate at the same position as, or neighbor to, index $i$ at level $h-1$. The restriction of constraints gives the lower bound constraints $\varphi_{h-1}$ at level $h-1$ in the form
\begin{equation}\label{constraint}
\begin{split}
\left( \varphi_{h-1}\right)_i & = (I_h^{h-1} x_{h,k})_i - \min \lbrace ( x_h)_j  |  j \in \mathcal I_{h,i} \rbrace
\end{split}
\end{equation}
This constraint is enforced to all iterates $x_{h-1,k}$ at the coarse level. Note that in FISTA at the target level, the constraint is enforced to $x_{h,k}$, whereas the transfer between levels is done at $y_{h,k}$.

\subsection{Decision on Recursive Search Direction}
At the beginning of each iteration $k$ the algorithm decides whether to compute a recursive search direction on the coarse level. This depends on the first-order optimality condition at the current iterate $y_{h,k}$ at the two levels, as well as the distance between the current iterate and the point $\tilde{y}_h$ at which the last recursive search direction was performed \cite{Wen,Parpas,Hovhannisyan}. In particular, a recursive search direction is used at $y_{h,k}$ if
\begin{eqnarray}  \label{conditions}
\providecommand{\norm}[1]{\left\lVert#1\right\rVert}
\begin{split}
&\left(\norm{I_h^{h-1} \nabla {F_\rho}_h(y_{h,k})} > \kappa \norm{\nabla {F_\rho}_h(y_{h,k})} \right) \hspace{0.2cm} \cap \hspace{0.2cm}  \\
&\left( \norm{y_{h,k}-\tilde{y}_h} > \vartheta \norm{\tilde{y}_h} \hspace{0.2cm}  \cup  \hspace{0.2cm} K_r=0  \hspace{0.2cm} \cup \hspace{0.2cm} K_d>q_d   \right),
\end{split}
\end{eqnarray}
where $ \providecommand{\norm}[1]{\left\lVert#1\right\rVert}
\kappa \in \left(0,\min(1,\min \norm{I_h^{h-1}})\right)$, $\vartheta \in (0,1)$, and $q_d \in \mathbb{N}$ are some predefined parameters, $K_d$ is the number of consecutive iterations with direct search direction, and $K_r$ is the number of all iterations already performed with recursive search direction \cite{Parpas,Hovhannisyan}. The first condition implies that a recursive search direction is not efficient if the first-order optimality condition is almost satisfied at the starting point of the coarse error correction, as this makes the minimization of the objective function at the coarse level ineffective. Furthermore, the second condition implies that a recursive search direction is not efficient if the current point $y_{h,k}$ is very close to the point $\tilde{y}_h$ since it gives a result that is similar to what was obtained on the last recursive search direction \cite{Wen}. This condition is ignored if the algorithm has already performed no step with a recursive search direction ($K_r=0$), or many consecutive steps with a direct search direction at the fine level, say greater than $q_d$.

\subsection{Outline of the MG algorithm}
In Algorithm \ref{FISTAMG} our MG variant of FISTA on two grids is outlined.
\begin{algorithm}
\caption {FISTA  in MG framework}
\begin{algorithmic}[1] \label{FISTAMG}
\STATE {\textbf{Iteration 0:} $y_{h,1}=x_{h,0}, \hspace{0.1cm} \theta_k=0, \hspace{0.1cm} K_d=0, \hspace{0.1cm} K_r=0$}
\STATE {\textbf{Iteration $\mathbf{k \geqslant 1}$:}}
\IF {$k>1 \quad \cap \quad \text{\eqref{conditions} holds}$}
\STATE   {Recursive search direction: $K_d=0,\hspace{0.1cm}  K_r=K_r+1$}
\STATE   {$x_{h-1,0}=I_h^{h-1} y_{h,k}$}
\STATE   {$v_{h-1}=I_h^{h-1} \nabla F_\rho(y_{h,k})-\nabla F_\rho (x_{h-1,0})$ }   
\STATE   {compute $\varphi_{h-1}$ by \eqref{constraint} }
\STATE   {compute $\phi_{h-1}$ by \eqref{fourteen}}
\STATE   {$x_{h-1,*}=\text{FISTA}(h-1,\phi_{h-1},x_{h-1,0},\varphi_{h-1})$}
\STATE   {$x_{h,k}=y_{h,k}+I_{h-1}^h\left(x_{h-1,*}-x_{h-1,0}\right)$}
\ELSE
\STATE   {Direct search direction: $K_d=K_d+1$}
\STATE   {$x_{h,k}=P_L(y_{h,k})$}
\ENDIF
\STATE   {update $\theta_{h,k}$}
\STATE   {$y_{h,k+1}  = x_{h,k}+ \theta_{h,k} \left(x_{h,k}-x_{h,k-1}\right)$}
\STATE   {\textbf{Output:}   $x_{h,*}$}
\end{algorithmic}
\end{algorithm}
Here, $*$ denotes the last iteration at each level. Since $\phi_{h-1}$ is smooth, at level $h-1$ the proximal map is reduced to a projection on the feasible set defined by $\varphi_{h-1}$. At each iteration with a recursive search direction, the termination of the algorithm at the coarse level is done whenever
\begin{align}
\frac{\phi_{h-1,k}-\phi_{h-1,k+1}}{\max \left(\phi_{h,k},\phi_{h,{k+1}}\right)}< \varepsilon_c \quad  \cup   \quad   *>q_c,
\end{align}
where $q_c$ denotes the maximum permitted number of iterations at the coarse level, and is applied to guarantee the efficiency of coarse error corrections.
Similarly, at the fine level the algorithm was terminated at iteration $k$ if
\begin{align}
\frac{F_{h,k}-F_{h,k+1}}{\max \left(F_{h,k}, F_{h,{k+1}}\right)}< \varepsilon_d.
\end{align}

\section{Numerical Results}
Numerical studies were performed to investigate the effectiveness of the proposed multi-grid strategy on performance of ISTA and FISTA
for iterative image reconstruction in PAT. To numerically solve
the three-coupled first order acoustic wave equations, which were described in Section \ref{aa}, the \textit{K-Wave} MATLAB toolbox was used \cite{Treeby-b}.
Additionally, to compute the gradient defined in (\ref{Six}) at each iteration at each level, the adjoint operator $H^*$ was computed based on the ``adjoint-then-discretize'' method (cf. Section \ref{A1}). The processor that was employed in this work is an Intel(R) Core(TM) i5-4570 CPU $@$ 3.20 GHz with a RAM of 8.00 GB and
a 64-bit operating system (Windows 7, Microsoft).
\subsection{2D PAT Simulation}   \label{s1}
A square grid with a size of $2.36 \times 2.36 \hspace{0.05cm}\text{cm}^2$ was created, which is made up of $472 \times 472$ grid points evenly
spaced with a separation distance of $5 \times 10^{-2}$mm in both $x$ and $y$ dimensions, supporting frequencies of up to $13.23$MHz. To measure the propagated wavefield, 200 point-wise pressure detectors were equidistantly placed along the left half of a circle having a radius of $11$mm so that $\pi$ radians of the circle were covered by the detectors. A PML having a thickness of $20$ grid points and a maximum attenuation coefficient of $2$ nepers per grid point was added to each side of the simulated grid in order to reduce spurious reflections at the boundaries \cite{Treeby-b}.

\begin{figure}\centering
{\subfigure[]{\includegraphics[scale=0.12]{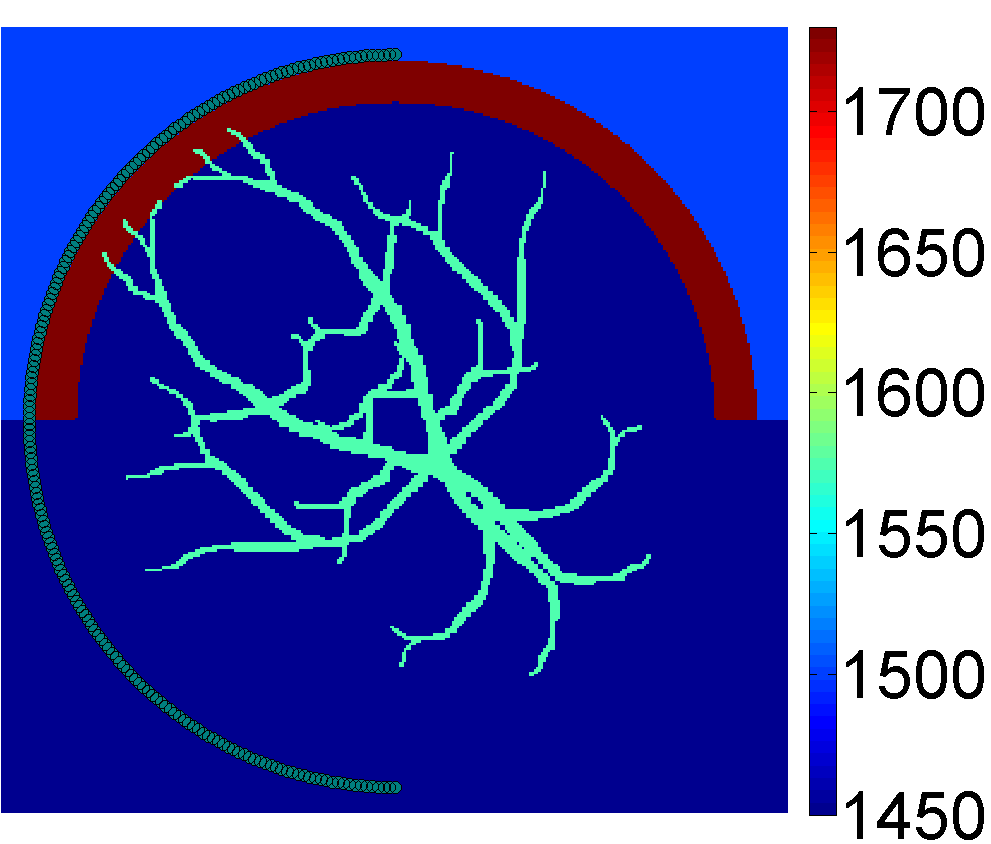}\label{F1a}}
\hspace{0.5cm}
\subfigure[]{\includegraphics[scale=0.12]{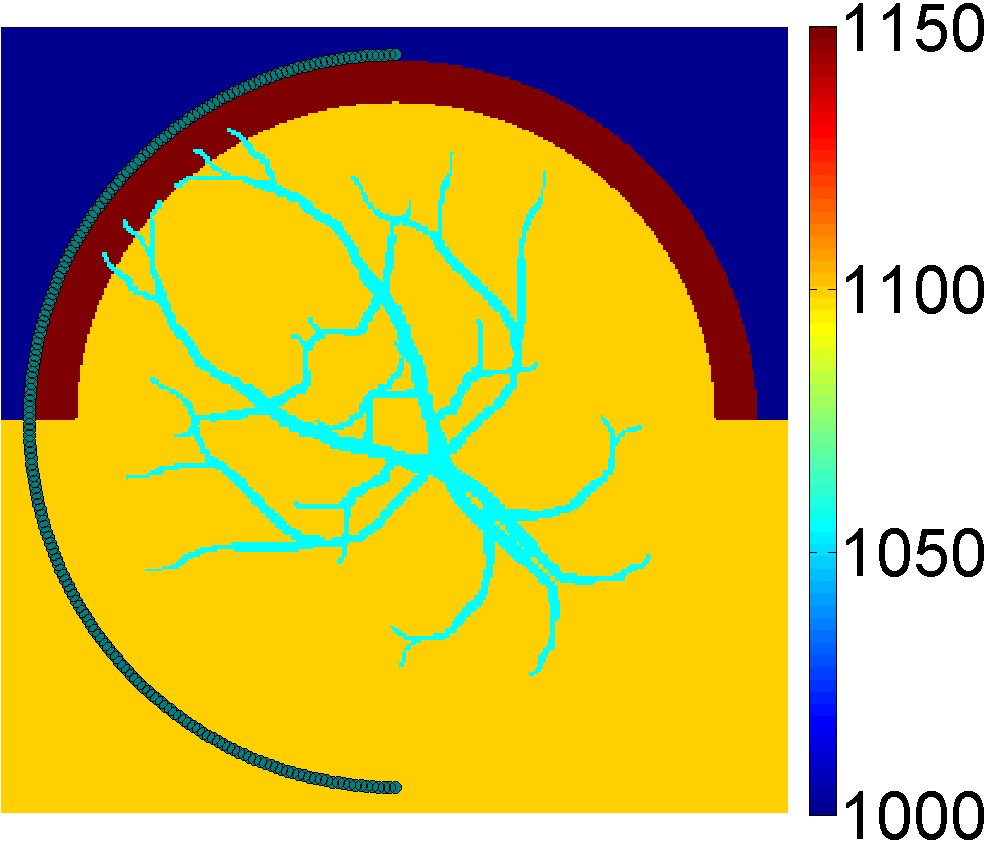}\label{F1b}}
\hspace{0.5cm}
\subfigure[]{\includegraphics[scale=0.12]{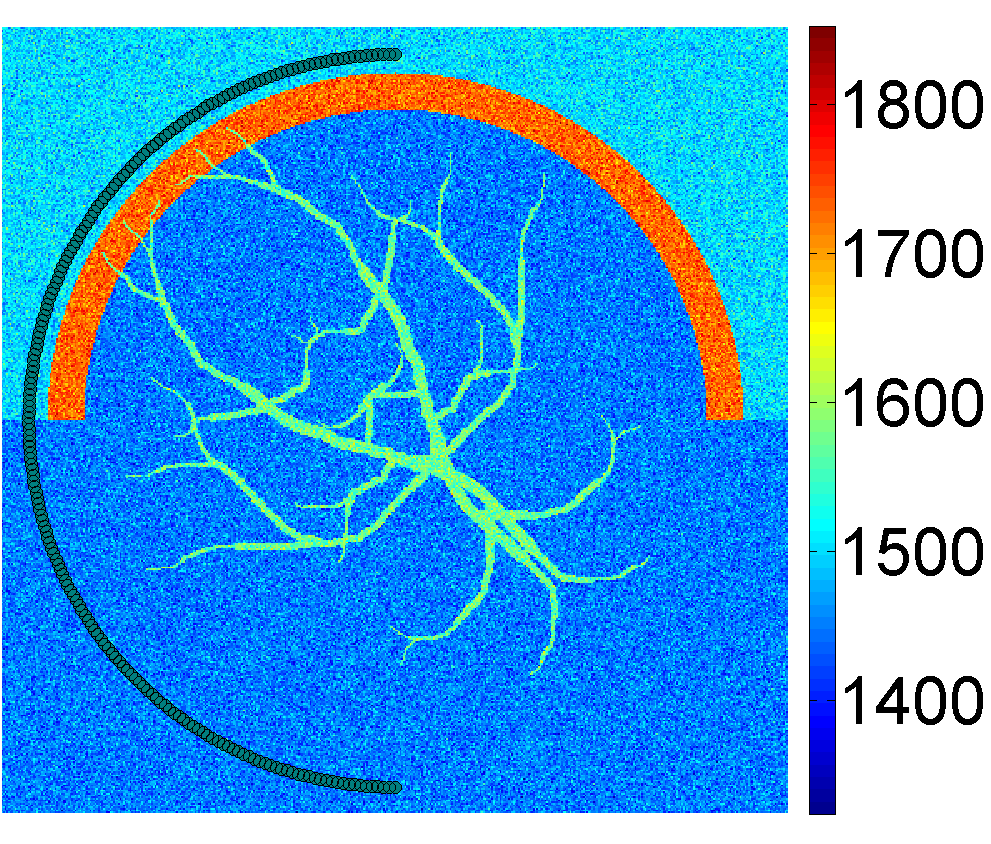}\label{F1c}}
\hspace{0.5cm}
\subfigure[]{\includegraphics[scale=0.12]{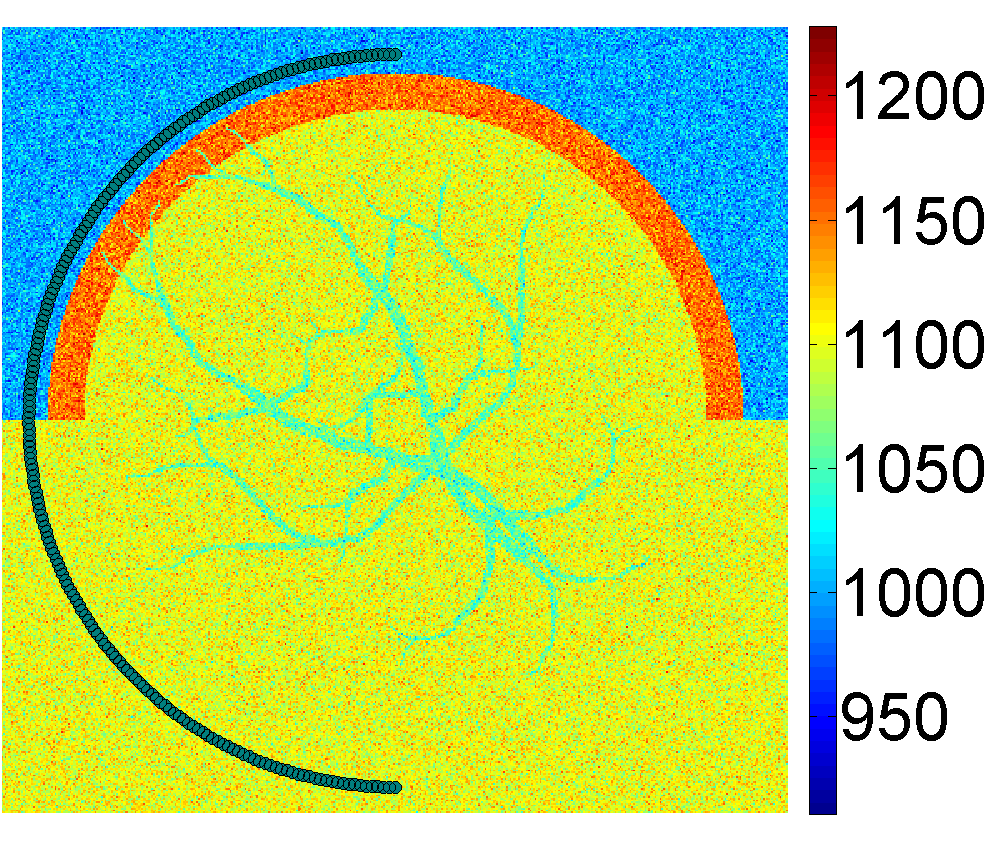}\label{F1d}}
}
\caption{Medium's properties for image reconstruction (a) sound speed (b) density, and data generation (c) sound speed (d) density.}
\end{figure}

\textbf{Medium's properties:} Figs. \ref{F1a} and \ref{F1b} show sound speed and density maps that were used for reconstruction, respectively. The sound speed and density for the inhomogeneity (vasculature) were set to those of blood, i.e., 1575 $\text{ms}^{-1}$ and 1055 $\text{kgm}^{-3}$, respectively, and the red color represents skin with a sound speed of 1730 $\text{ms}^{-1}$ and a density of 1150 $\text{kgm}^{-3}$. The background inside the detection surface represents fat tissue with a sound speed of 1450 $\text{ms}^{-1}$ and a density of 950 $\text{kgm}^{-3}$, and a sound speed of 1500 $\text{ms}^{-1}$ and a density of 1000$\text{kgm}^{-3}$ were considered for region outside the detection surface to represent water. These maps were inspired by acoustic properties in tissues given in \cite{Azhari}. Note that acoustic properties in realistic tissues are often smoother than the simulated maps, and do not have sharp interfaces. However, these sharp maps were provided in order to make a challenge for coarse error correction in the MG method.

Furthermore, the absorption coefficient was set to 0.75 $\text{dB} \hspace{0.1cm} \text{MHz}^{-y} \hspace{0.1cm} \text{cm}^{-1}$ for the whole medium, except the area that represents water, where it was set to $2 \times 10^{-3} \hspace{0.1cm} \text{dB} \hspace{0.1cm} \text{MHz}^{-y} \hspace{0.1cm} \text{cm}^{-1}$. The attenuation power law exponent was set to $1.5$ for the entire medium.

Since the exact maps are not readily available for reconstruction, data were generated from a more realistic phantom by contaminating the maps with a $35$dB AWGN, as well as shifting the ``water-skin'' and ``skin-soft tissue'' interfaces towards the centre of the detection surface by $2\%$ of radius of the circle.
Figs. \ref{F1c} and \ref{F1d} show distributions of sound speed and density that were used for data generation, respectively. To mitigate errors arising from aliasing, for all forward and adjoint models, the acoustic properties were smoothed by the k-wave toolbox \cite{Treeby-b}.

\begin{figure}\centering
{\subfigure[]{\includegraphics[scale=0.12]{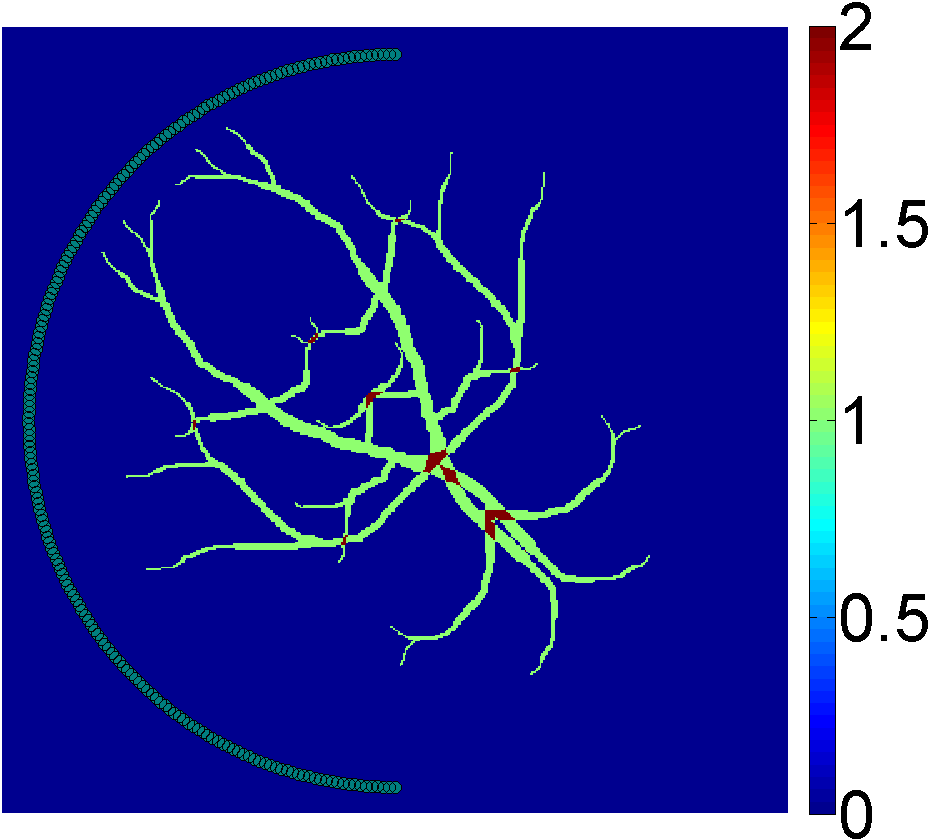}\label{F2a}}
\hspace{0.5cm}
\subfigure[]{\includegraphics[scale=0.12]{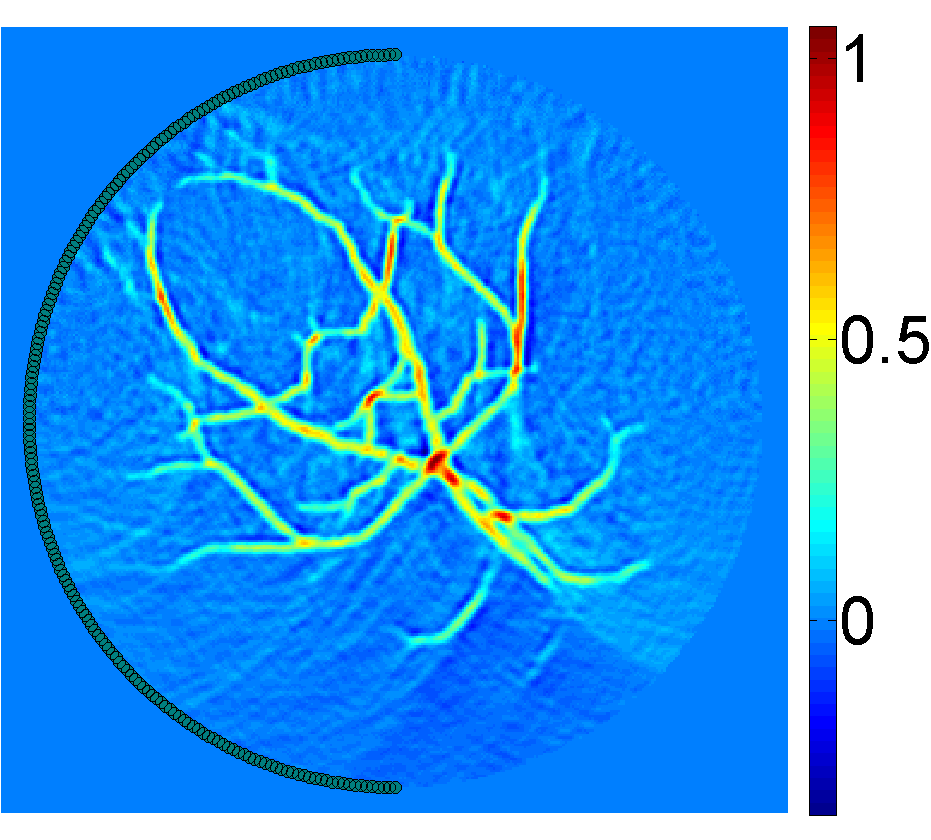}\label{F2b}}
}
\caption{2D phantom. (a) Initial pressure distribution (b) image reconstructed by TR.}
\end{figure}

The phantom was created so that it simulates the pressure distribution of vessels with a maximal amplitude of 2. Fig. \ref{F2a} displays the simulated phantom. To avoid spurious oscillations in the computed pressure field, high frequencies of the initial pressure distribution for each forward implementation were filtered by a self-adjoint smoothing operator. This operator was then included in the adjoint (see \cite{Arridge}).

The computed time-dependent pressure field arriving at the sensors was then sampled evenly in $2655$ time steps. The interpolation of pressure field to sensors was performed by the well-known linear method. The generated data was then contaminated with a $30$ dB AWGN. In order to avoid inverse crime, the reconstruction was applied to a grid made up of $328\times 328$ grid points, which supports a maximal frequency of $10.07$ MHz. The PML at each side of the grid was proportionally reduced to $16$ grid points. This grid will further be used as the ``target grid'' for our proposed multi-grid algorithm.

\textit{Iterative methods}: The iterative reconstruction was performed by TV-regularized ISTA, i.e., Algorithm \ref{ISTA}.
The step size $\alpha_k$ was chosen to be $2/L_f$ \cite{Beck-b} computed by the ``power iteration'' method, similar to \cite{Arridge,Arridge-b}. The regularization parameter was heuristically set to $\lambda=1 \times 10^{-2}$.

The MG variant of ISTA, Algorithm \ref{FISTAMG} with $\text{FISTA}$ replaced by $\text{ISTA}$, was then employed to reconstruct images on two grids having sizes $328\times 328$ and $164\times 164$. The algorithm was implemented by $\kappa=1/4$, $\vartheta=10^{-1}$, $q_d=3$, $q_c=8$, $\varepsilon_d=10^{-3}$ and $\varepsilon_c=10^{-2}$.
The coarse model supports a maximal frequency half the fine grid, i.e. $5.038$ MHz. For iterations at which a recursive search direction was computed, the TV function was smoothed by $\rho=1 \times 10^{-2}$ as in \eqref{Thirty-two}. At the coarse level, the number of grid points associated with the PML was halved so that the thickness of the PML was the same as the target level. The sequences $x^k$ computed by the algorithms were measured by the following parameters.

\textit{Relative Error (RE)}: This is defined at iteration $k$ as
\begin{align}
\providecommand{\norm}[1]{\left\lVert#1\right\rVert}
RE(x^k)=\frac{\norm{p_\text{sol}-p_\text{exact}}}{\norm{p_\text{exact}}} \times 100,
\end{align}
where $p_\text{sol}$ stands for the sequence $x^k$ interpolated back to the forward grid, and $p_\text{exact}$ denotes the simulated phantom.

\textit{ Norm of Residual (RES)}: This is defined at iteration $k$ as
\begin{align}\providecommand{\norm}[1]{\left\lVert#1\right\rVert}
\text{RES}(x^k)=\norm{Ax^k-\hat{p}}.
\end{align}

\textit{Objective function ($F$)}: This is defined at iteration $k$ as a discretized variant of \eqref{One} on the inverse grid. It should be noted that the efficiency of any optimization algorithm, including the MG algorithm we are examining, should be evaluated using the objective function.

\begin{figure}\centering
{\subfigure[]{\includegraphics[scale=0.08]{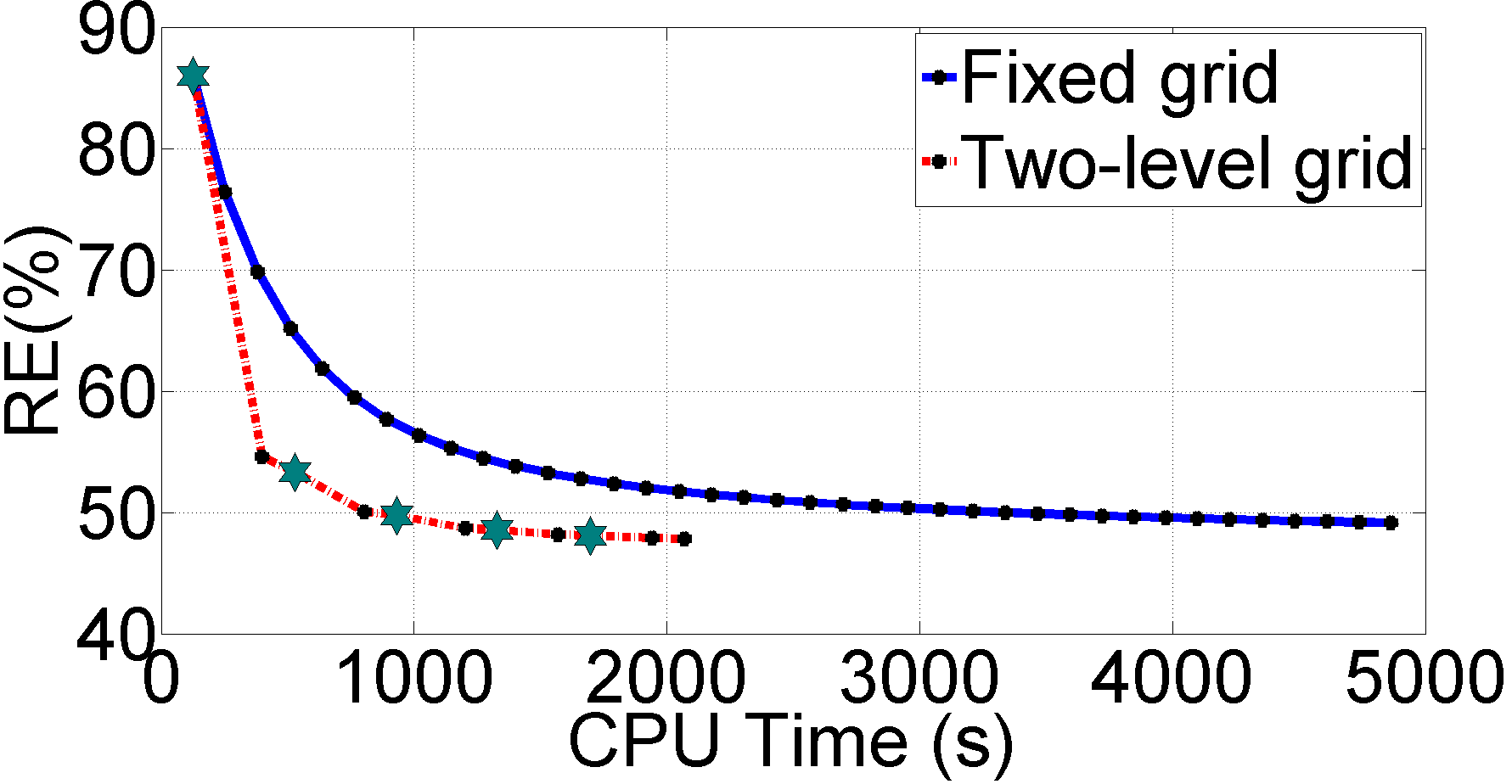}\label{F3a}}
\subfigure[]{\includegraphics[scale=0.08]{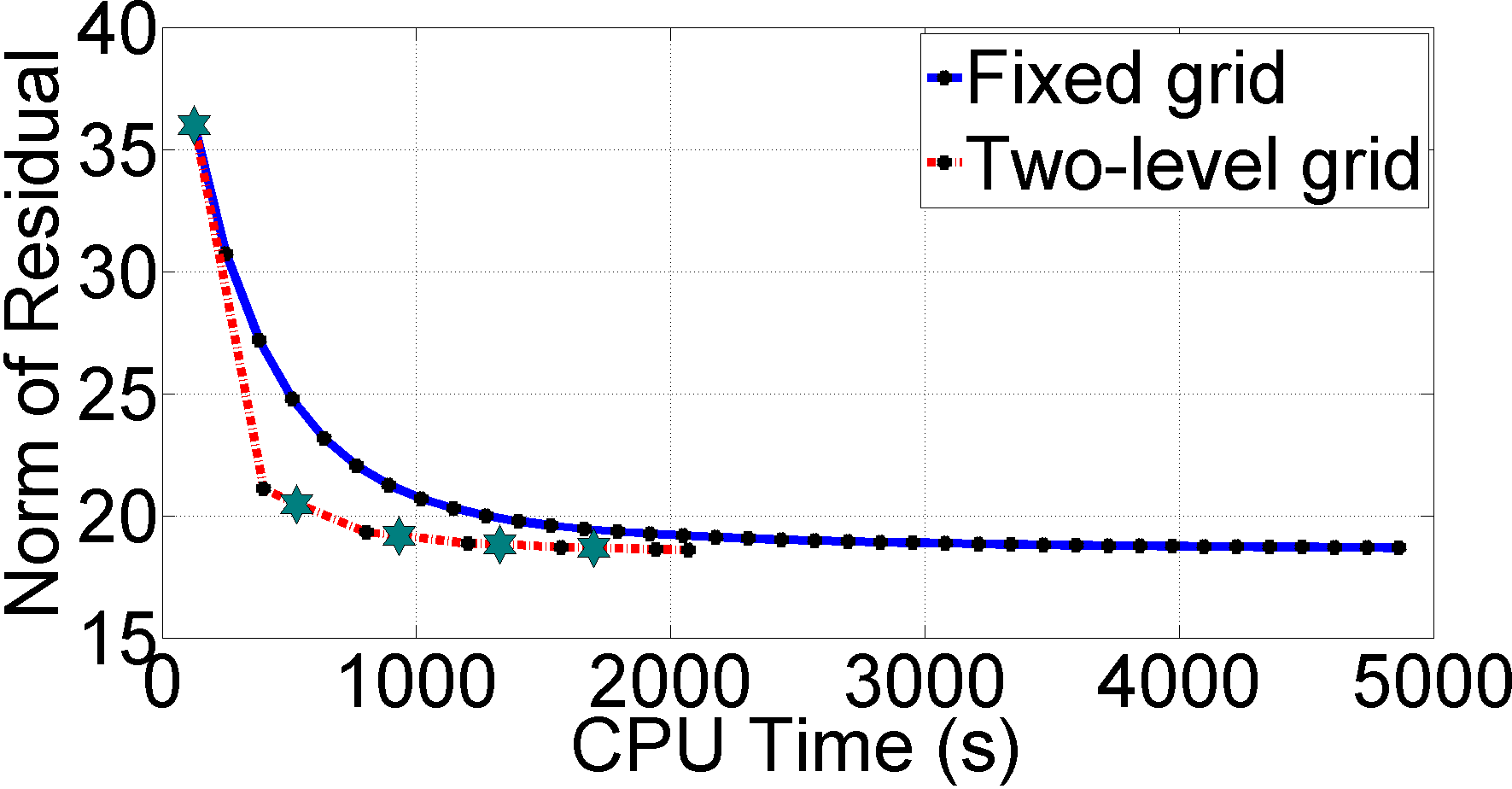}\label{F3b}}\\
\subfigure[]{\includegraphics[scale=0.08]{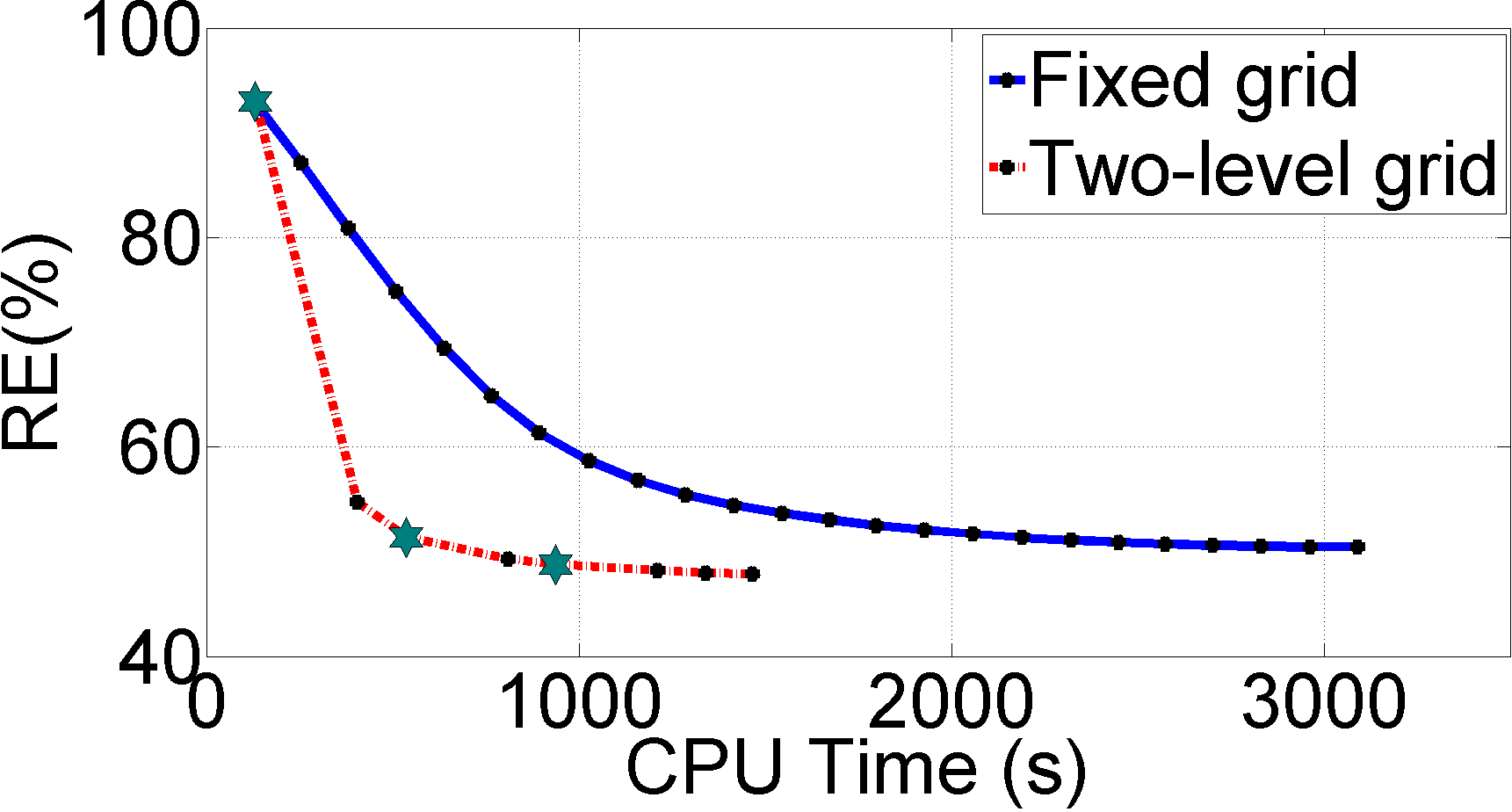}\label{F3c}}
\subfigure[]{\includegraphics[scale=0.08]{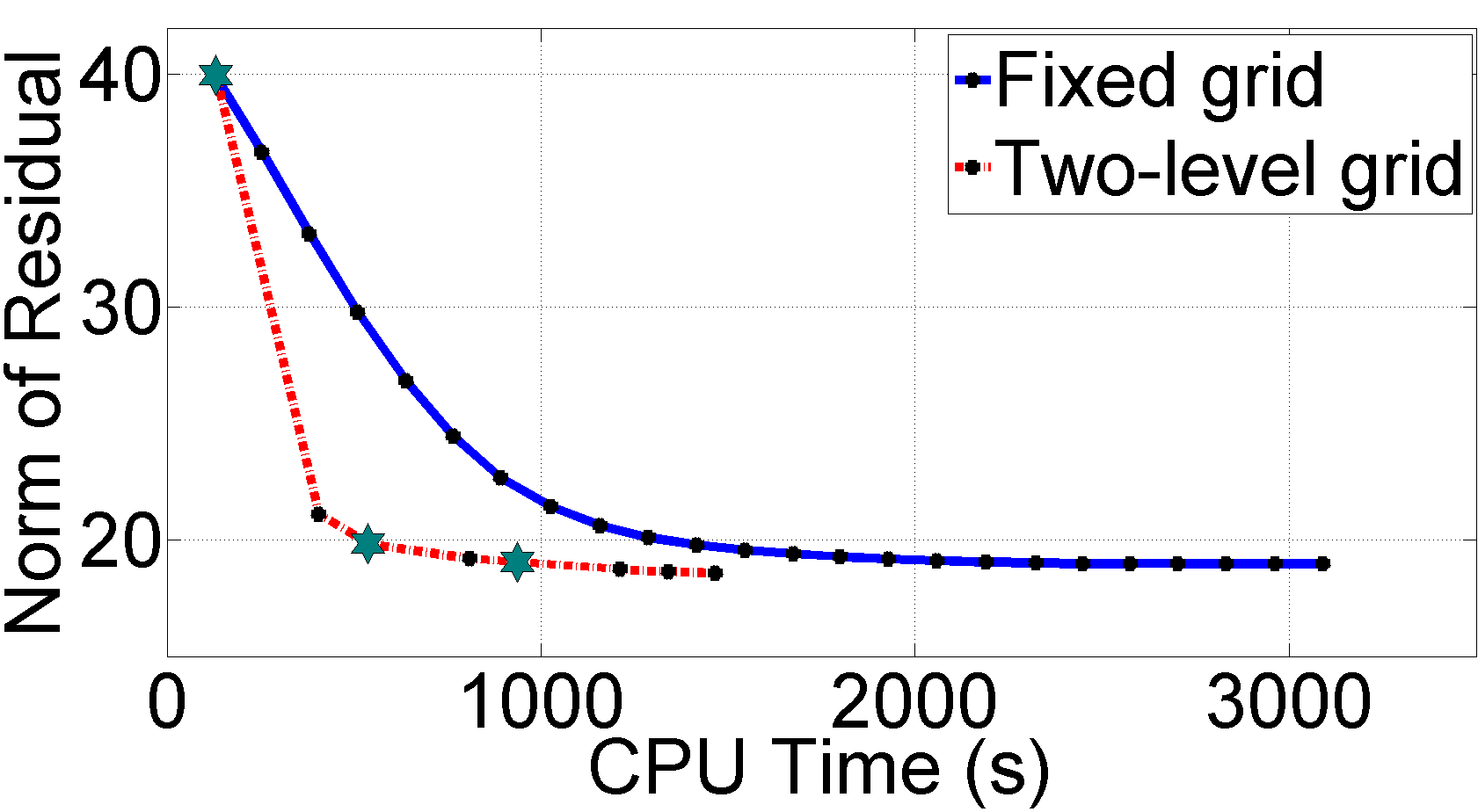}\label{F3d}}
}
\caption{Evaluation of 2D images reconstructed on fixed grid (blue) and two-level grid (red). ISTA: (a) RE (b) RES, and FISTA: (c) RE (d) RES.}
\end{figure}

Fig. \ref{F2b} shows the image reconstructed by TR, which has an $RE$ of $68.95\%$. Fig. \ref{F3a} shows RE of the images reconstructed by ISTA. To make a fair comparison between the competing algorithms, the image parameters were plotted versus CPU time although the iterations are also shown as black dots. For the MG algorithm, the iterations at which a recursive search direction was used are designated by hexagrams. As seen in this figure, the fixed-grid algorithm was terminated after $4.87 \times 10^3$s (38 iterations), and finally reconstructed an image having an RE of $49.14\%$, whereas the MG algorithm produced an image having an RE of $48.68\%$ at $1.20 \times 10^3$s (6 iterations). The MG algorithm was finally terminated after $2.07 \times 10^3$s (11 iterations), and provided an image having an RE of $47.81\%$. Fig. \ref{F3b} shows RES in the same way as RE.

Subsequently, FISTA and its MG variant were implemented on the same grids as and with the same parameters as ISTA and MG ISTA, respectively, except that the step size was chosen to be $1/L_f$ \cite{Beck,Beck-b}. Figs. \ref{F3c} and \ref{F3d} display RE and RES of sequences provided by FISTA, respectively. As shown in Fig. \ref{F3c}, FISTA reconstructed a final image having an RE of $50.44\%$ after $3.09 \times 10^3$s on a fixed grid (24 iterations), while the MG variant of FISTA provided an image having an RE of $49.31\%$ at $8.09 \times 10^2$s ($4$ iterations). The MG algorithm finally reached an RE of $47.87\%$ at $1.46 \times 10^3$s ($8$ iterations).

\begin{figure}\centering
{\subfigure[]{\includegraphics[scale=0.08]{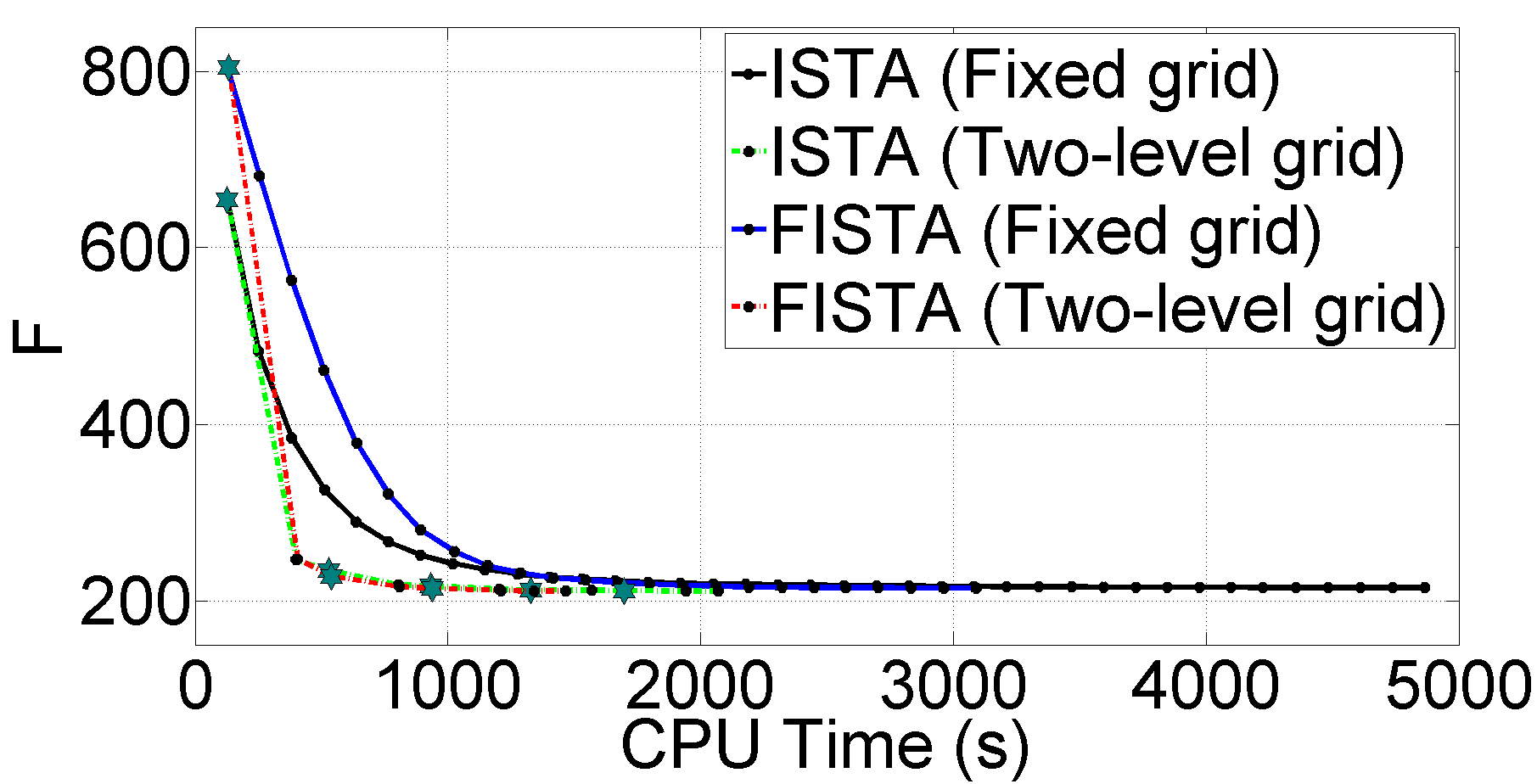}\label{F4a}}
\subfigure[]{\includegraphics[scale=0.08]{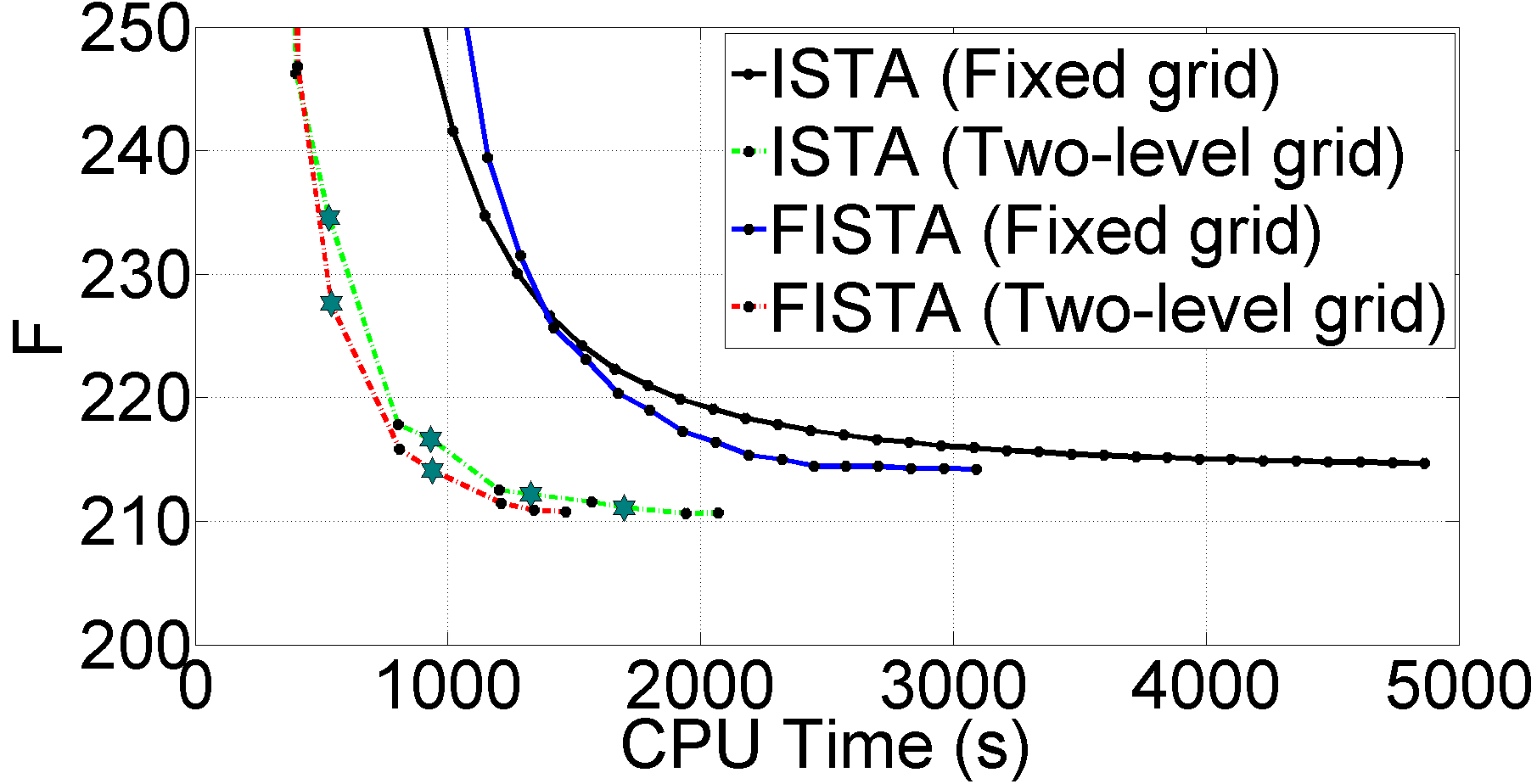}\label{F4b}}
}
\caption{Convergence of ISTA on a fixed grid (black), ISTA on a two-level grid (green), FISTA on a fixed grid (blue) and FISTA on a two-level grid (red) for the 2D phantom. (a) objective function (b) larger view.}
\end{figure}

The objective function values corresponding to sequences computed by ISTA and FISTA are shown in Fig. \ref{F4a}. This figure has been shown in a larger view around the optimum in Fig. \ref{F4b}. As seen in these figures, FISTA converged more slowly than ISTA at early iterations, but it was faster than ISTA around the optimum on both fixed grid and two-level grids. ISTA has finally reached an objective function having a value of $2.15 \times 10^2$ at $4.87 \times 10^3$s on a fixed grid, whereas it
has reached an $F$ of $2.12 \times 10^2$ at $1.20 \times 10^3$s on a two-level grid. This indicates that the MG variant of ISTA was
four times faster than the fixed-grid ISTA. The MG variant of ISTA has finally reached an $F$ of $2.10 \times 10^2$ at $2.07 \times 10^3$s. Applying FISTA on a fixed grid, the algorithm was terminated at $3.09 \times 10^3$s with an $F$ of $2.14 \times 10^2$, while FISTA on a two-level grid has reached almost the same value at $9.37\times 10^2$s (5 iterations). This implies that FISTA on a two-level grid was almost three times faster than on a fixed grid. MG FISTA was finally terminated at $1.46 \times 10^3$s with an $F$ of $2.10 \times 10^2$.

Some of the reconstructed images have been shown in Fig. 5. The color scale of each figure was set independently. It is clear that the MG variant of FISTA reconstructed more accurate images in less times than the fixed-grid algorithm.
\begin{figure}\centering
{\subfigure[]{\includegraphics[scale=0.12]{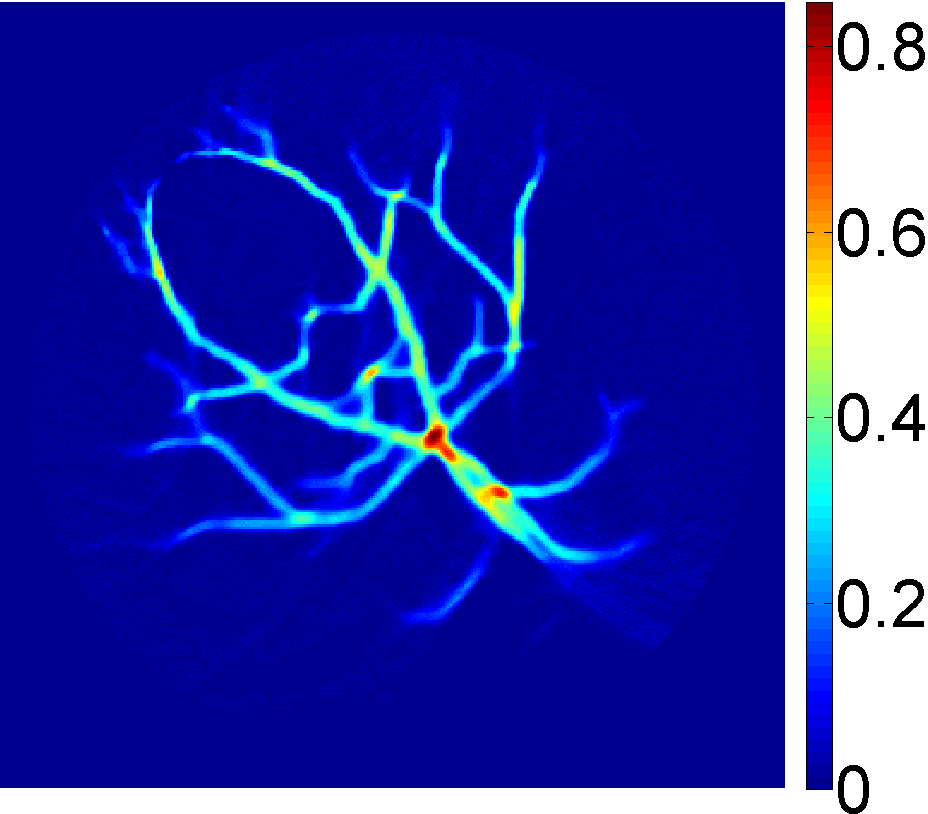}\label{F5a}}
\hspace{0.5cm}
\subfigure[]{\includegraphics[scale=0.12]{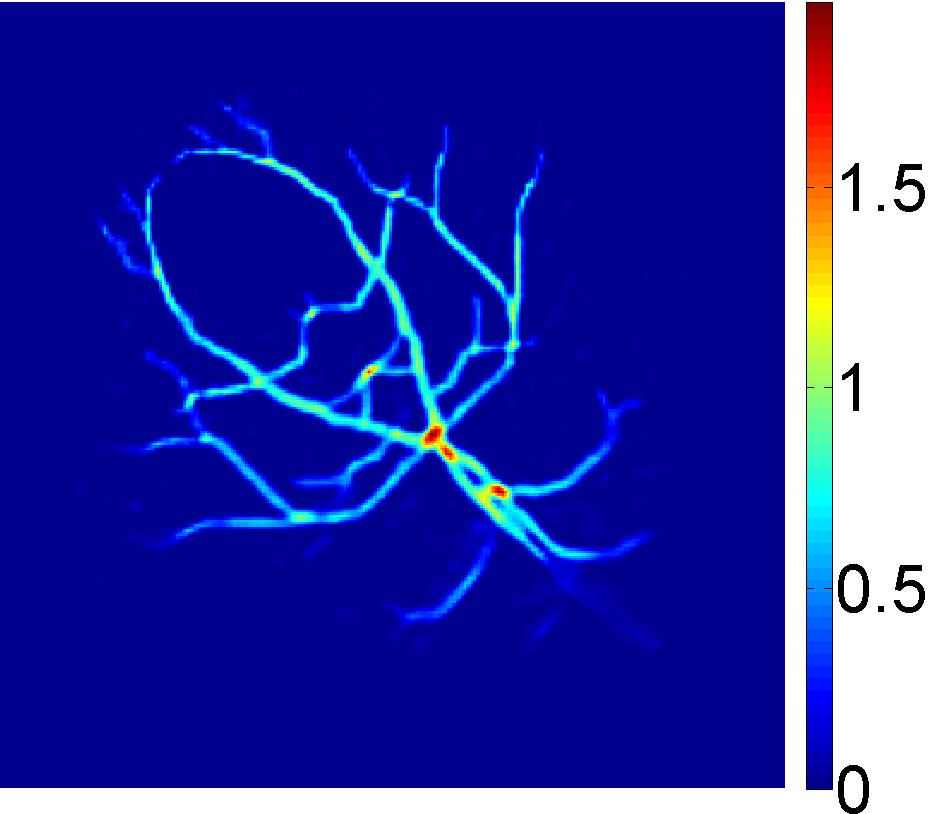}\label{Fb}}
\hspace{0.5cm}
\subfigure[]{\includegraphics[scale=0.12]{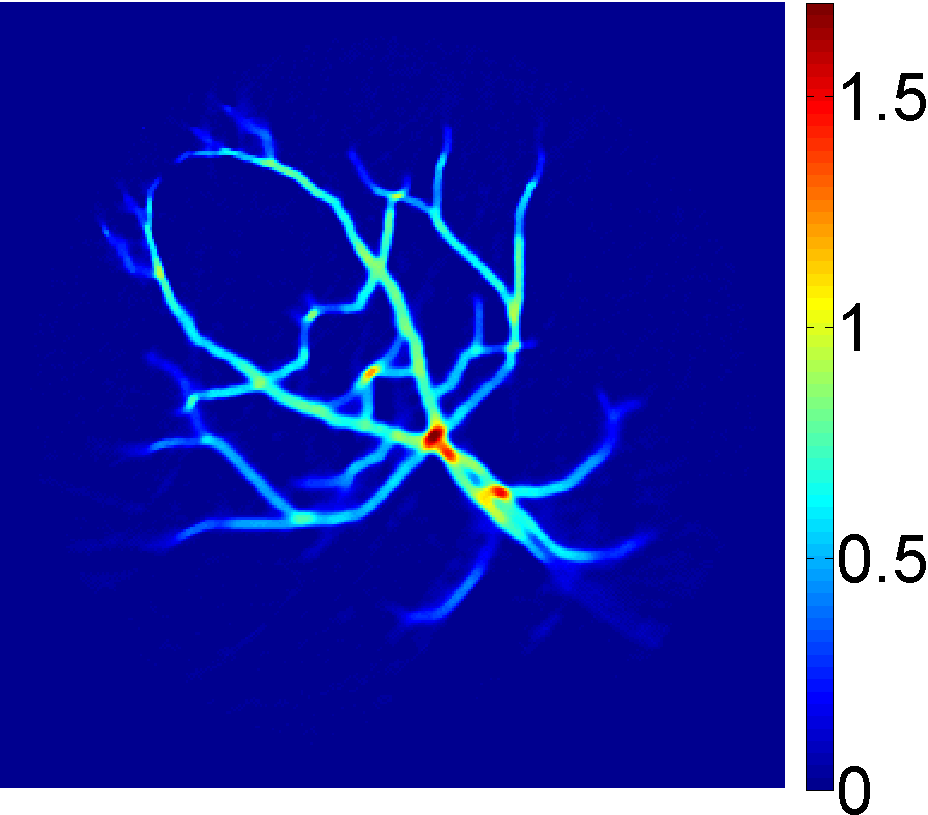}\label{F5c}}
\hspace{0.5cm}
\subfigure[]{\includegraphics[scale=0.12]{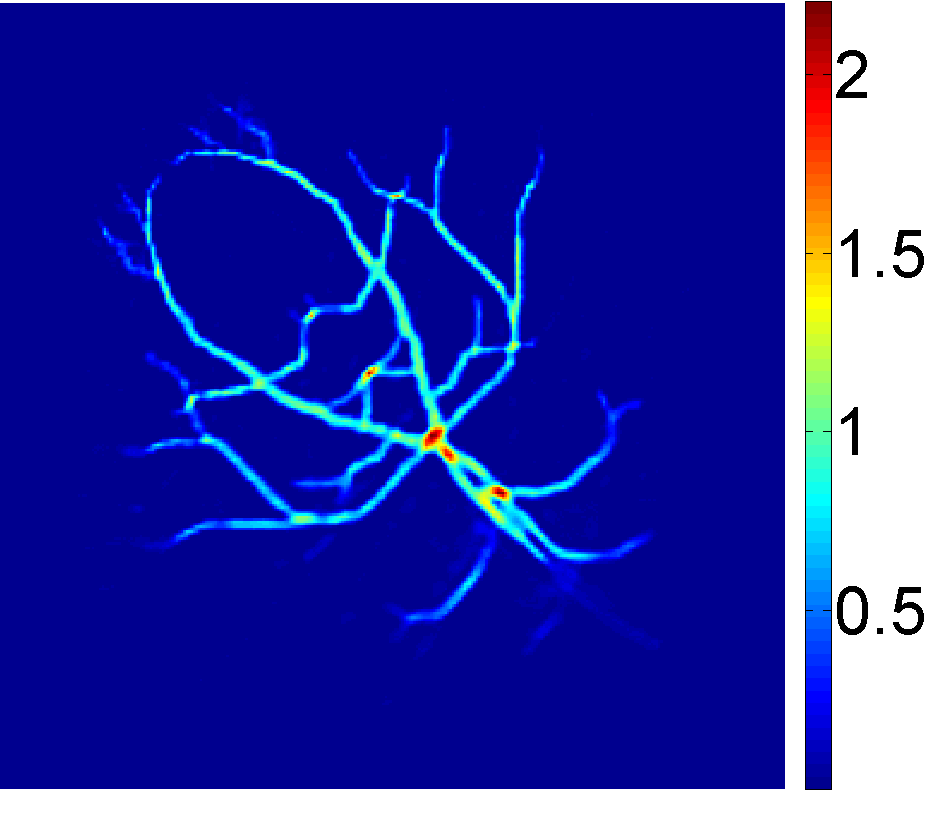}\label{F5d}}
}
\caption{2D images reconstructed by FISTA. (a) Iteration 4 on a fixed grid ($5.09 \times 10^2$s) (b) Iteration 2 on a two-level grid ($4.05 \times 10^2$s) (c) Iteration 8 on a fixed grid ($1.03 \times 10^3$s) (d) Iteration 4 on a two-level grid ($8.09 \times 10^2$s).}
\end{figure}

\subsection{3D PAT Simulation}
\begin{figure}\centering
{\subfigure[]{\includegraphics[scale=0.08]{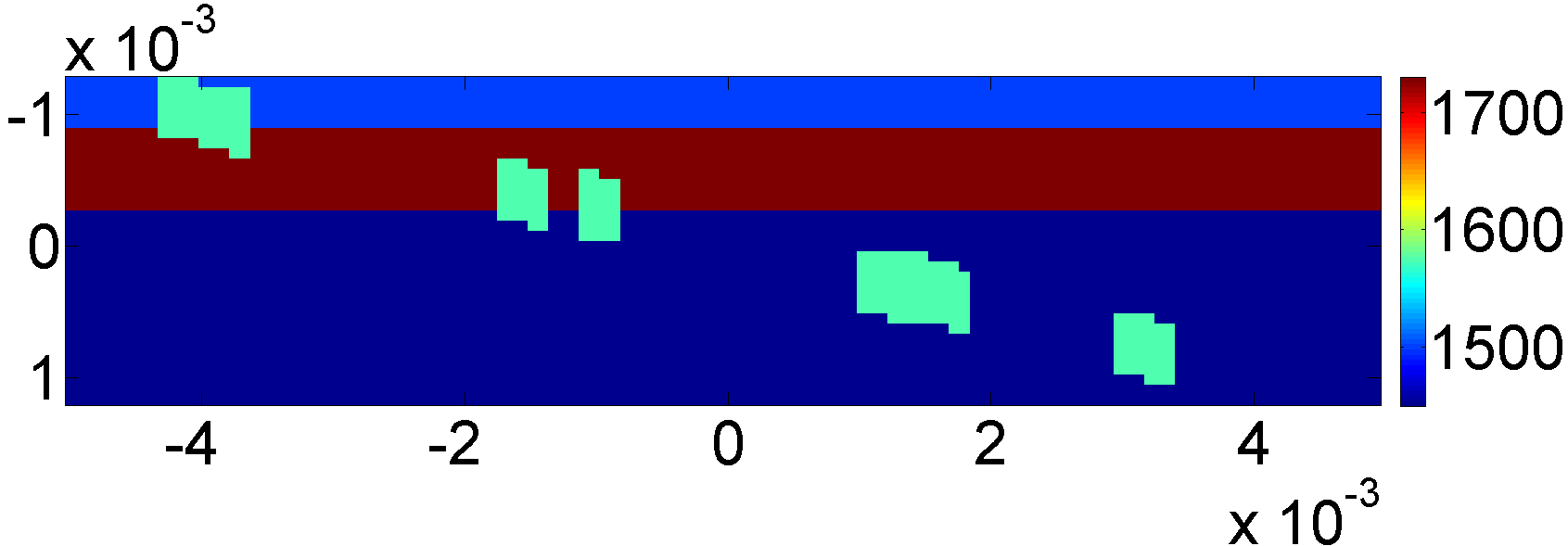}\label{F6a}}
\hspace{0.5cm}
\subfigure []{\includegraphics[scale=0.08]{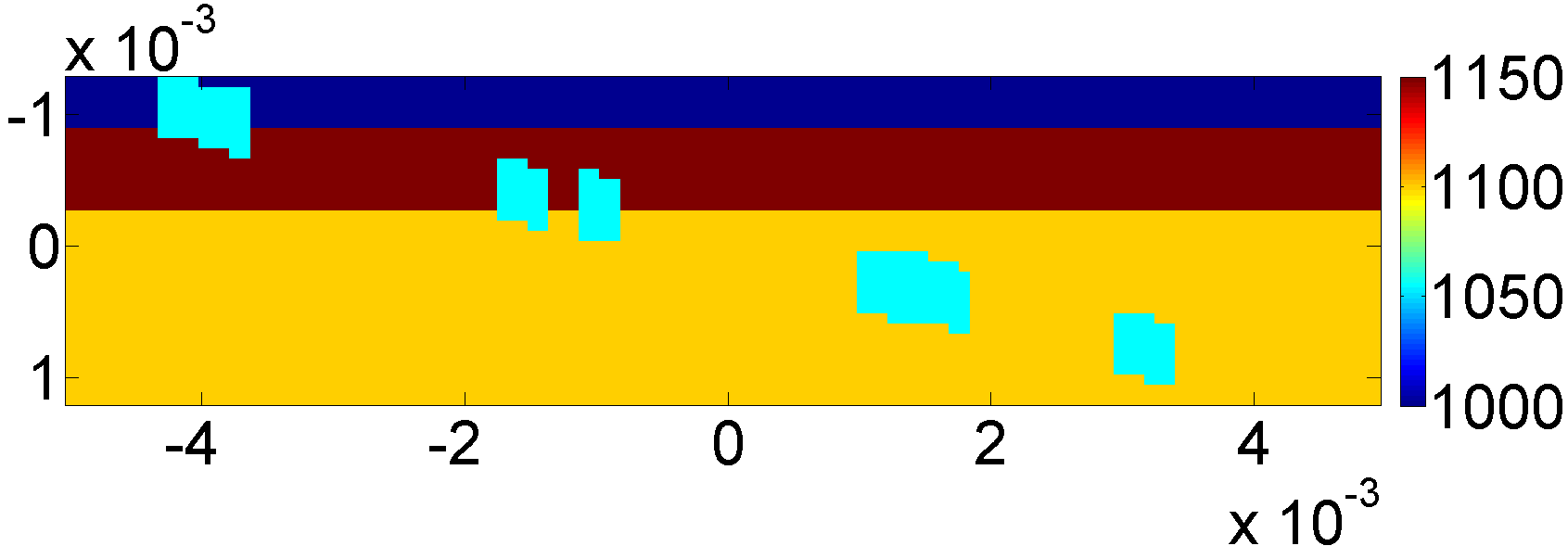}\label{F6b}}
\hspace{0.5cm}
\subfigure[]{\includegraphics[scale=0.08]{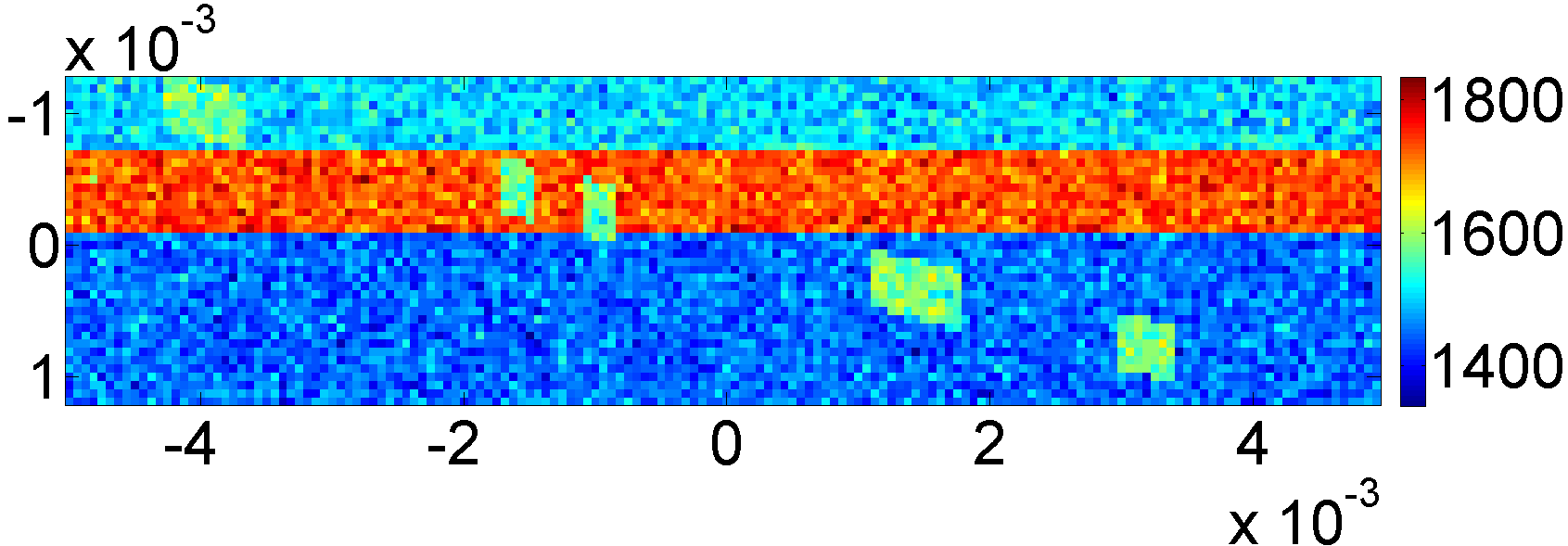}\label{F6c}}
\hspace{0.5cm}
\subfigure []{\includegraphics[scale=0.08]{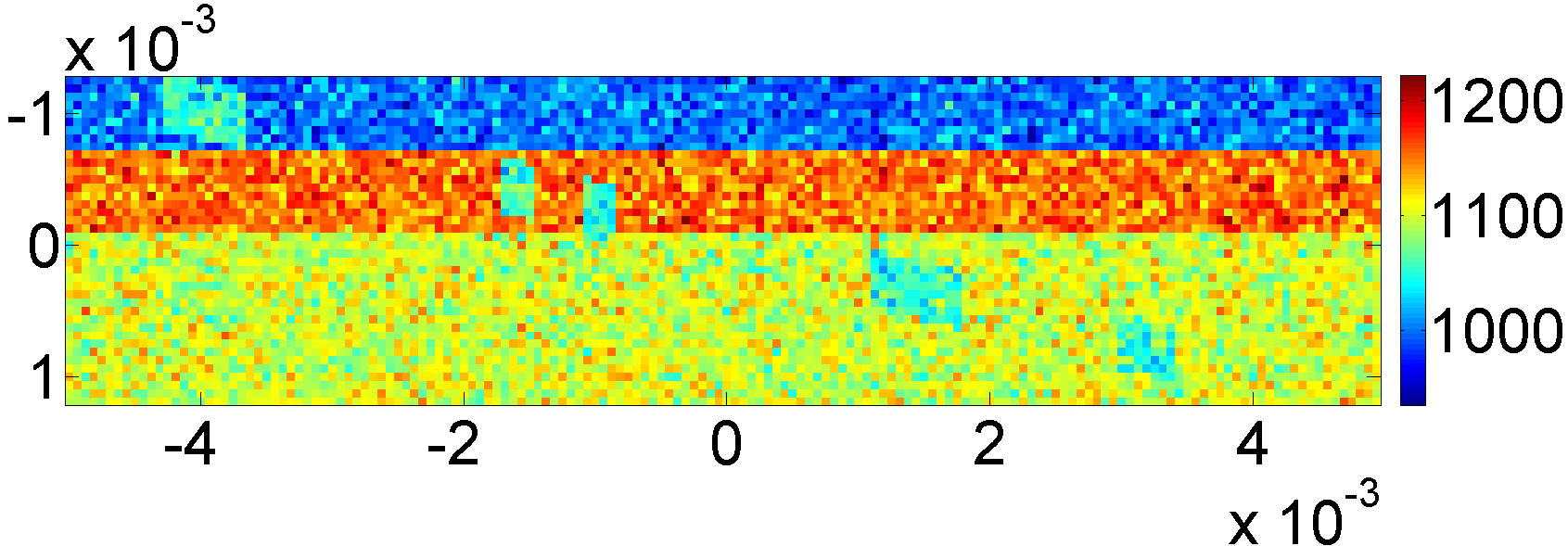}\label{F6d}}
}
\caption{Medium's properties for image reconstruction (a) sound speed (b) density, and data generation (c) sound speed (d) density.}
\end{figure}

A 3D grid with a size of $1 \times 1 \times 0.25  \hspace{0.05cm} \text{cm}^3$ was created, made up of $160 \times 160\times 40$ grid points with a spatial spacing of $6.25 \times 10^{-2}$ mm, supporting a maximal frequency of $10.36$ MHz in all axes. To measure the propagated wavefield, $36 \times 36$ point-wise pressure detectors were placed on the top surface of the grid (see \cite{Arridge-b}). A PML was added to the grid in the same way as in Section \ref{s1}.
The sound speed and density maps were simulated inhomogeneous, as shown in Figs. \ref{F6a} and \ref{F6b}, respectively. From the top to bottom, the layers represent properties of water, skin and soft tissue with sound speed and density values given in Section \ref{s1} \cite{Azhari}. The sound speed and density of vasculature were set the same as the 2D phantom as well.

To avoid inverse crime, for data generation these maps have been contaminated with a $35$dB AWGN, and the ``water-skin'' and ``skin-soft tissue'' interfaces were shifted to the bottom by $3$ grid points ($18.75 \times 10^{-2}$mm), as shown in Figs. \ref{F6c} and \ref{F6d}. To mitigate aliasing artifacts, for all forward and adjoint models, medium's properties were smoothed by the k-wave toolbox \cite{Treeby-b}. The absorption coefficient and power law exponent of tissues were set the same as those in Section \ref{s1}.

\begin{figure}\centering
{\subfigure[]{\includegraphics[scale=0.13]{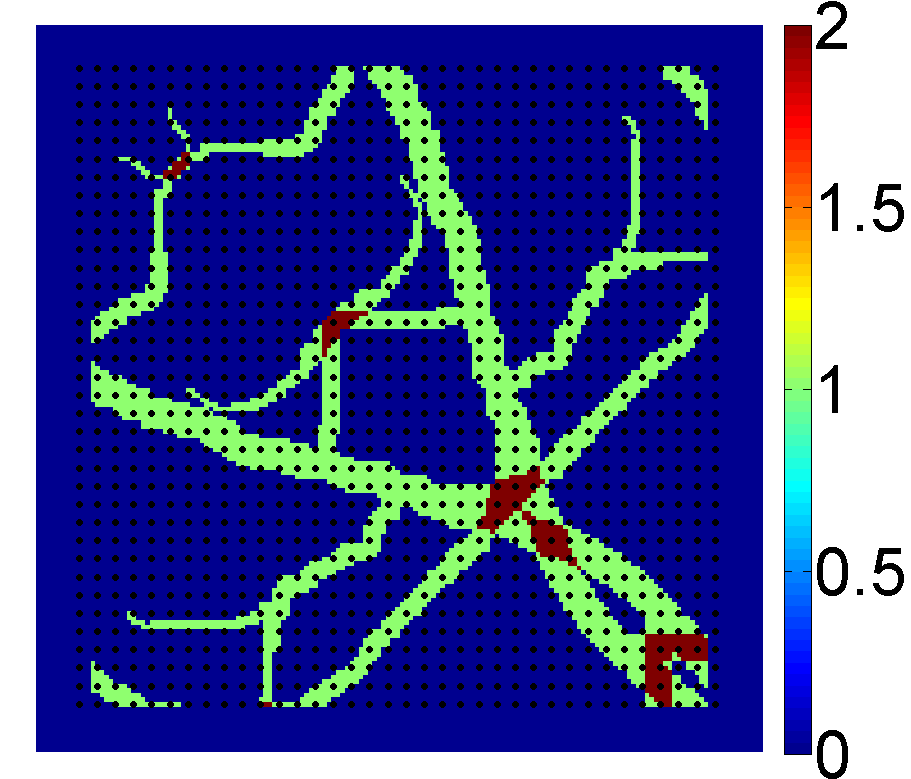}\label{F7a}}
\hspace{0.5cm}
\subfigure[]{\includegraphics[scale=0.13]{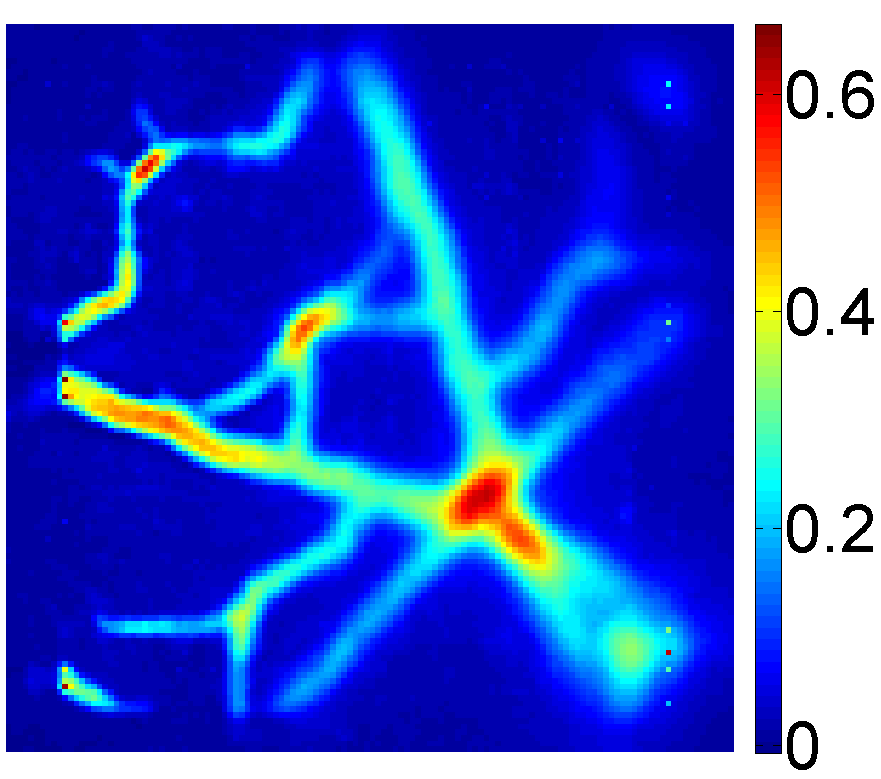}\label{F7b}}
}
\caption{3D phantom. (a) Initial pressure distribution (b) image reconstructed by TR.}
\end{figure}

The phantom that was already created for the 2D scenario was now placed obliquely in the plane $z=y/4$. Fig. \ref{F7a} displays the simulated phantom from a top view, and the sensors are shown by black dots. The computed pressure field was evenly sampled in $914$ time steps, and was linearly interpolated to the sensors. A $30$dB AWGN was then incorporated to the generated data.

The reconstruction was applied to a grid made up of $128 \times 128\times 32$ grid points, supporting a maximal frequency of $9.28$ MHz.
Fig. \ref{F7b} displays the image reconstructed by TR. This image has an RE of $87.98\%$.
In our study, 3D visualizations were done by Maximum Intensity Projection (MIP) method (see \cite{Huang-c,Arridge-b}).

\textit{Iterative methods:} The iterative reconstruction was performed by ISTA and FISTA. The step sizes were chosen by power iteration method similarly to the algorithms in Section \ref{s1}. The regularization parameter was heuristically chosen to be $\lambda=1 \times 10^{-2}$. The MG algorithms were then implemented to reconstruct images on two levels having sizes $128\times 128\times 32$ and $64 \times 64\times 16$. These algorithms were implemented by the same parameters as in 2D scenario, except that $\kappa$ was set to $1/8$, i.e., size of the coarse grid relative to the fine grid, and also the smoothing parameter $\rho$ was set to $3 \times 10^{-2}$.

\begin{figure}\centering
{\subfigure[]{\includegraphics[scale=0.08]{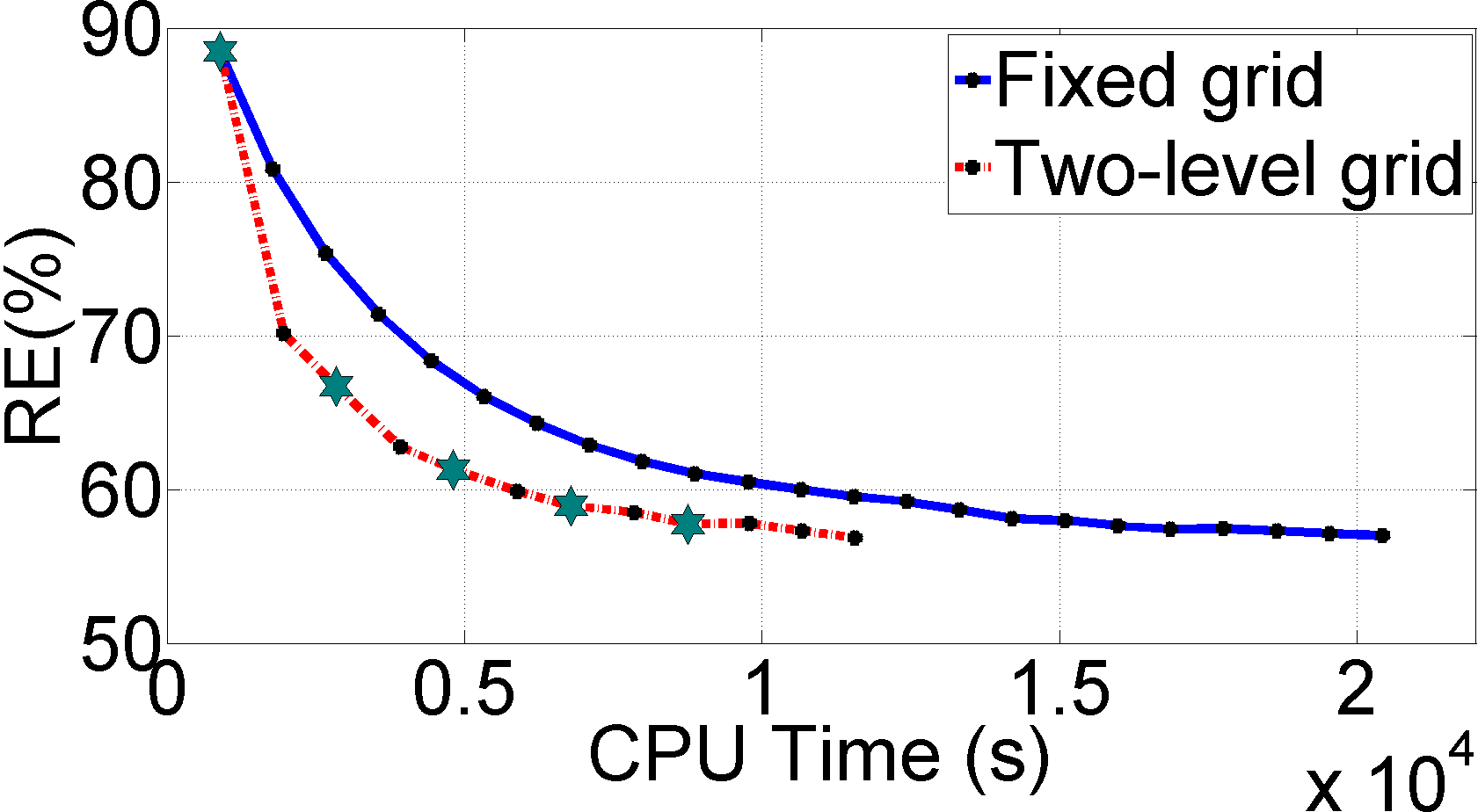}\label{F8a}}
\subfigure[]{\includegraphics[scale=0.08]{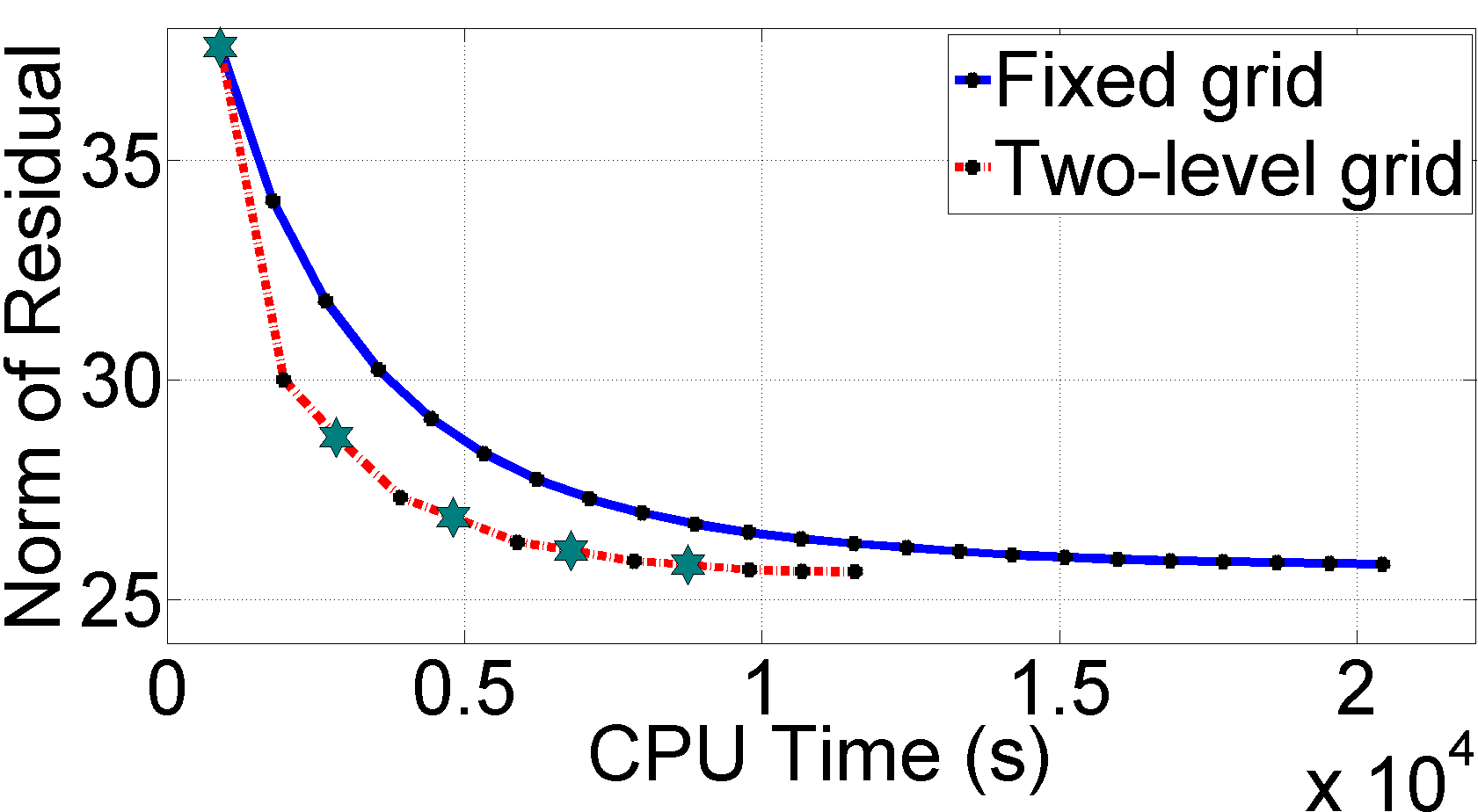}\label{F8b}}\\
\subfigure[]{\includegraphics[scale=0.08]{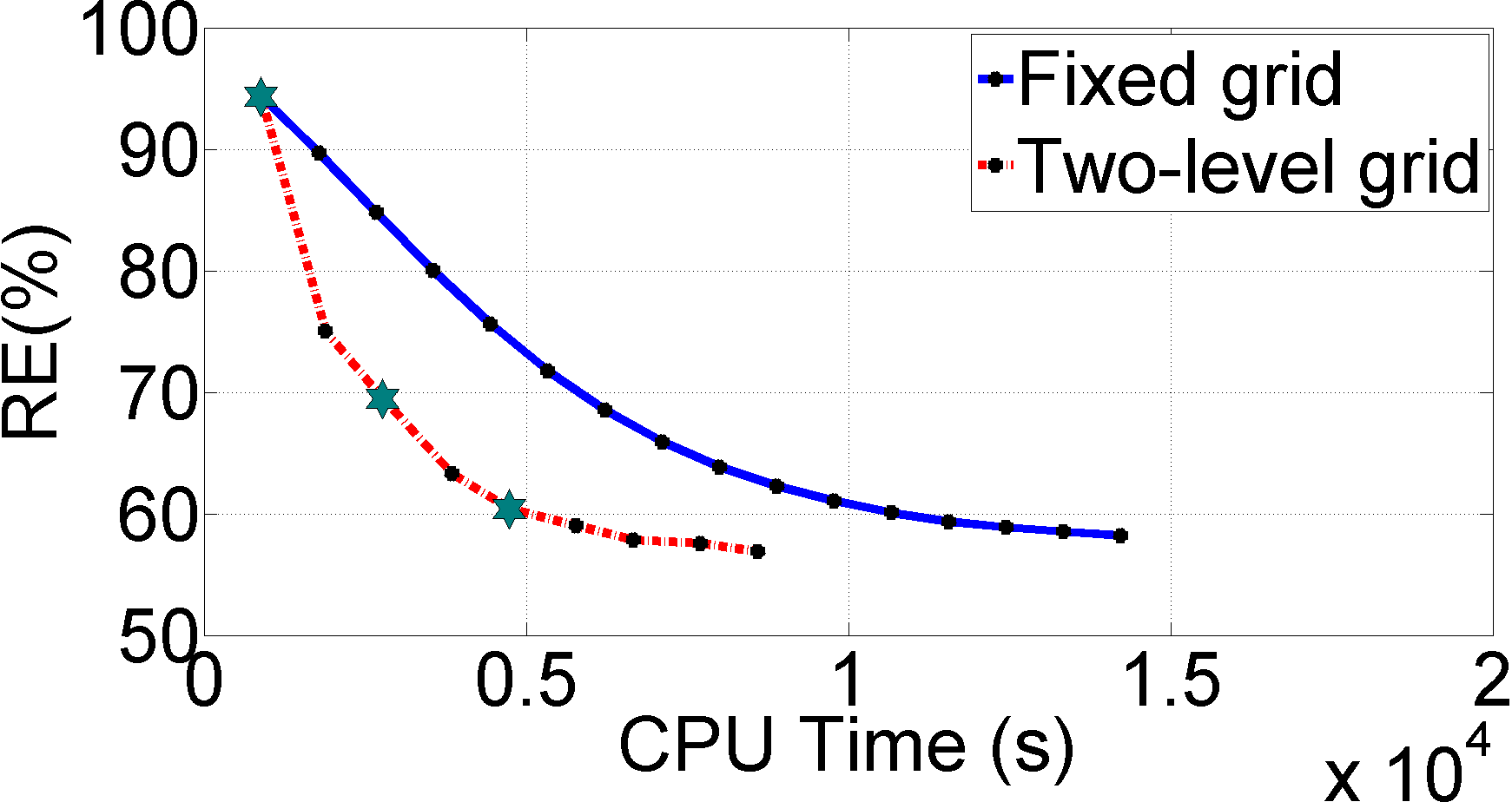}\label{F8c}}
\subfigure[]{\includegraphics[scale=0.08]{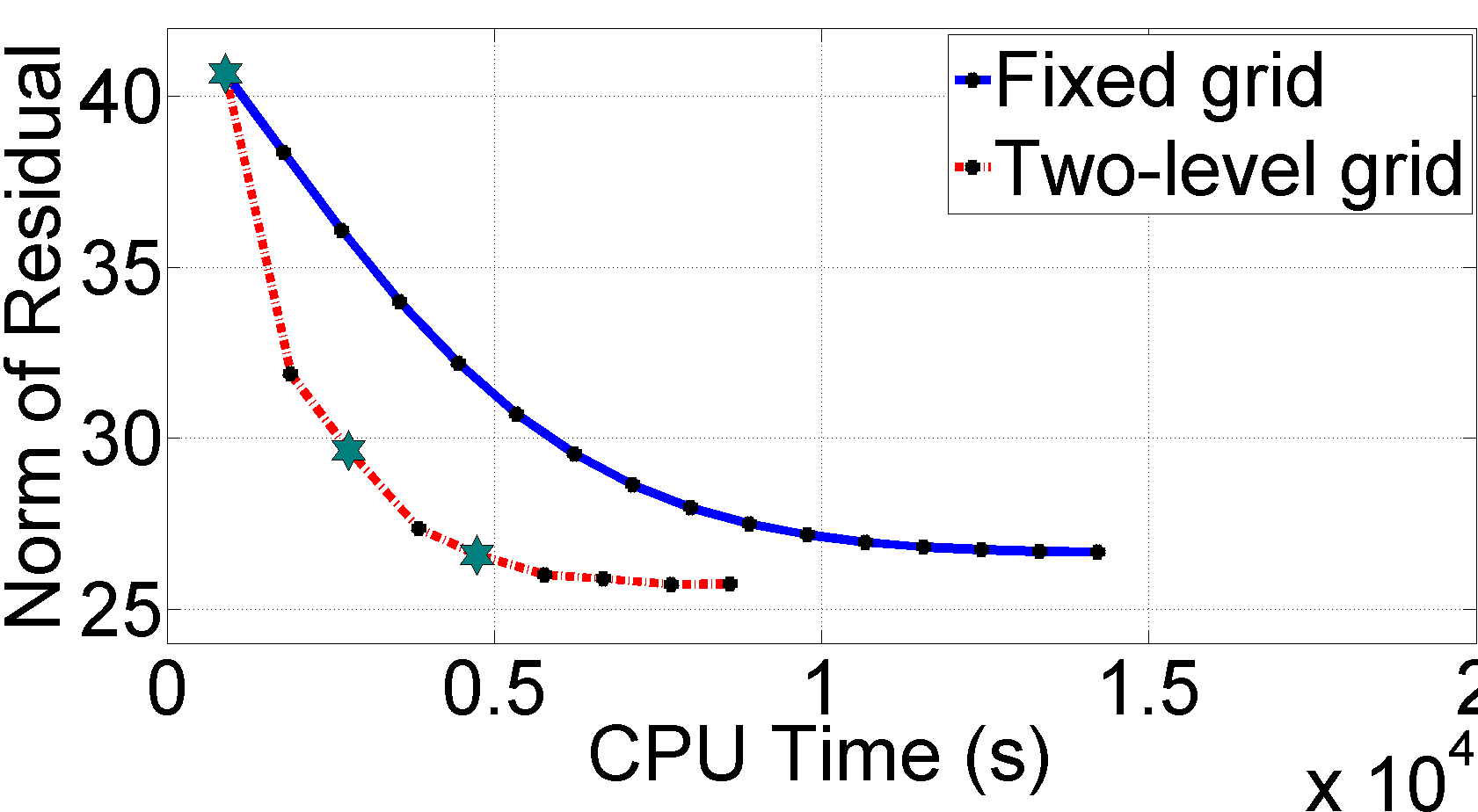}\label{F8d}}
\caption{Evaluation of 3D images reconstructed on fixed grid (blue) and two-level grid (red). ISTA: (a) RE (b) RES, and FISTA: (c) RE (d) RES.}}
\end{figure}

Figs. \ref{F8a} and \ref{F8b} respectively, show RE and RES of the images reconstructed by ISTA versus CPU time in the same way as in Section \ref{s1}.
Figs. \ref{F8c} and \ref{F8d} display RE and RES of sequences provided by FISTA, respectively. As seen in these figures, both ISTA and FISTA exhibited a better performance on a two-level grid than on a fixed grid.

\begin{figure}\centering
\includegraphics[scale=0.10]{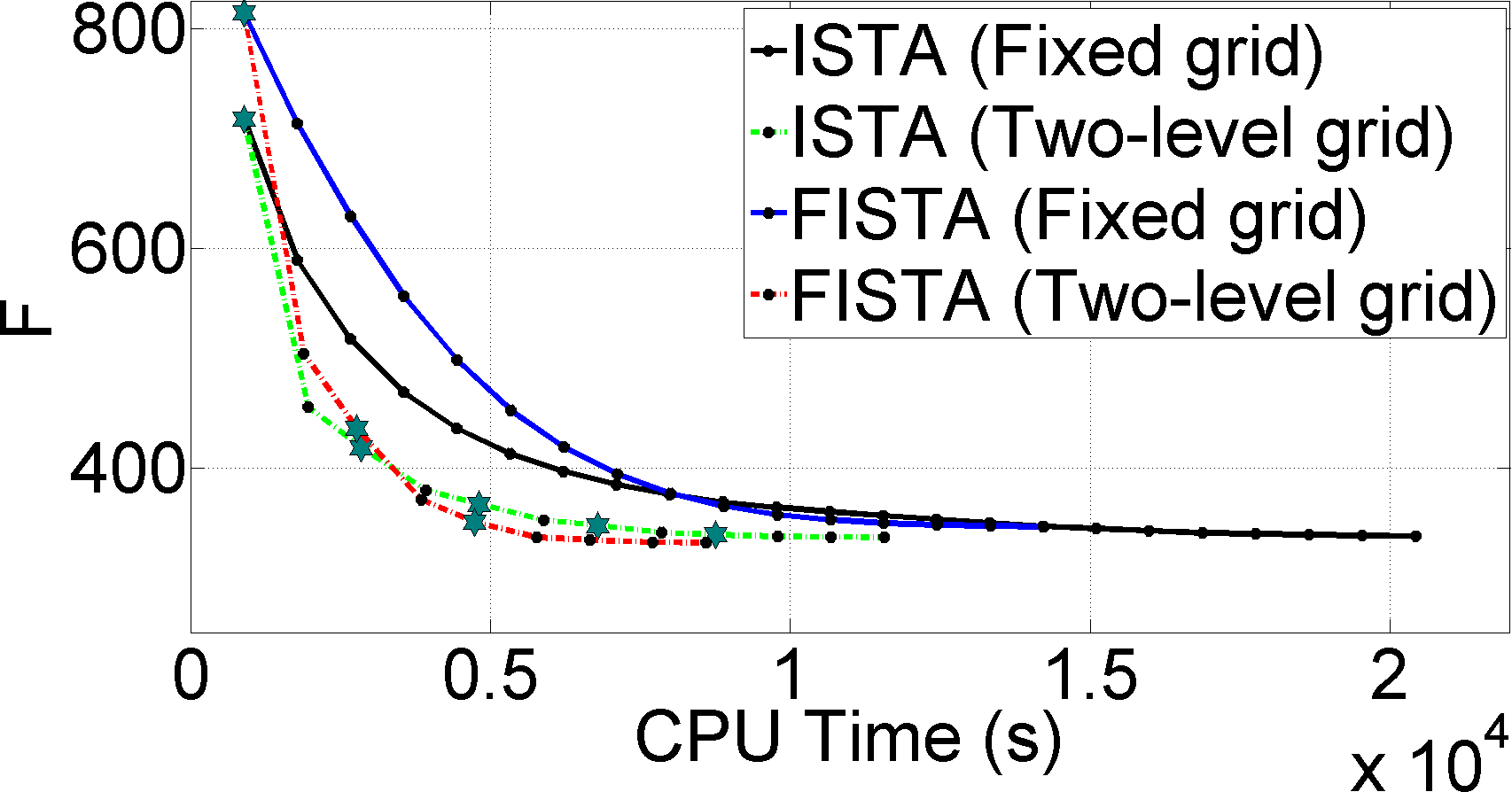}
\caption{Objective function values computed by ISTA on a fixed grid (black), ISTA on a two-level grid (green), FISTA on a fixed grid (blue), and FISTA on a two-level grid (red) for the 3D phantom.}
\label{F9}
\end{figure}

Fig. \ref{F9} shows $F$ values versus CPU time for all the used algorithms.
As shown in this figure, the $F$ value obtained by MG ISTA after $8.74 \times 10^3$s ($9$ iterations) was less than the optimal value obtained by fixed-grid ISTA at $2.05 \times 10^4$s (23 iterations). This indicates that MG ISTA was $2.35$ times faster than ISTA on a fixed grid.

Applying FISTA, the $F$ value computed by the MG algorithm after $4.71 \times 10^3$s ($5$ iterations) was less than the optimal value computed by the fixed-grid algorithm at $1.42 \times 10^4$s ($16$ iterations). This indicates that MG FISTA was almost three times faster than FISTA on a fixed grid.

Furthermore, a comparison between ISTA and FISTA on a fixed grid indicates that the convergence of ISTA was faster than FISTA at early iterations because of using a step size greater than FISTA, but FISTA has finally provided almost the same optimal $F$ as ISTA in a less time. Applying the MG algorithms, MG ISTA was faster than MG FISTA at early iterations, but MG FISTA converged better than MG ISTA around the optimal point, and reached a lower $F$ than MG ISTA in less time. MG ISTA finally reached a value of $3.36 \times 10^2$ at $1.15 \times 10^4$s (12 iterations), while MG FISTA was terminated at $8.59 \times 10^3$s (9 iterations) with an $F$ of $3.32 \times 10^3$.

\begin{figure}\centering
{\subfigure[]{\includegraphics[scale=0.13]{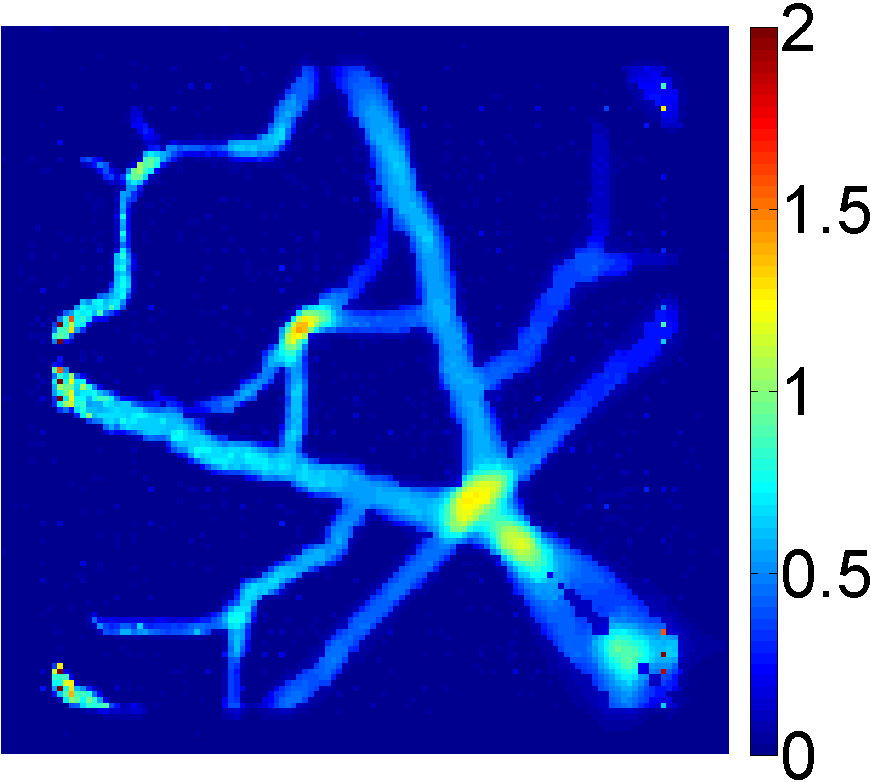}\label{F10a}}
\hspace{0.5cm}
\subfigure[]{\includegraphics[scale=0.13]{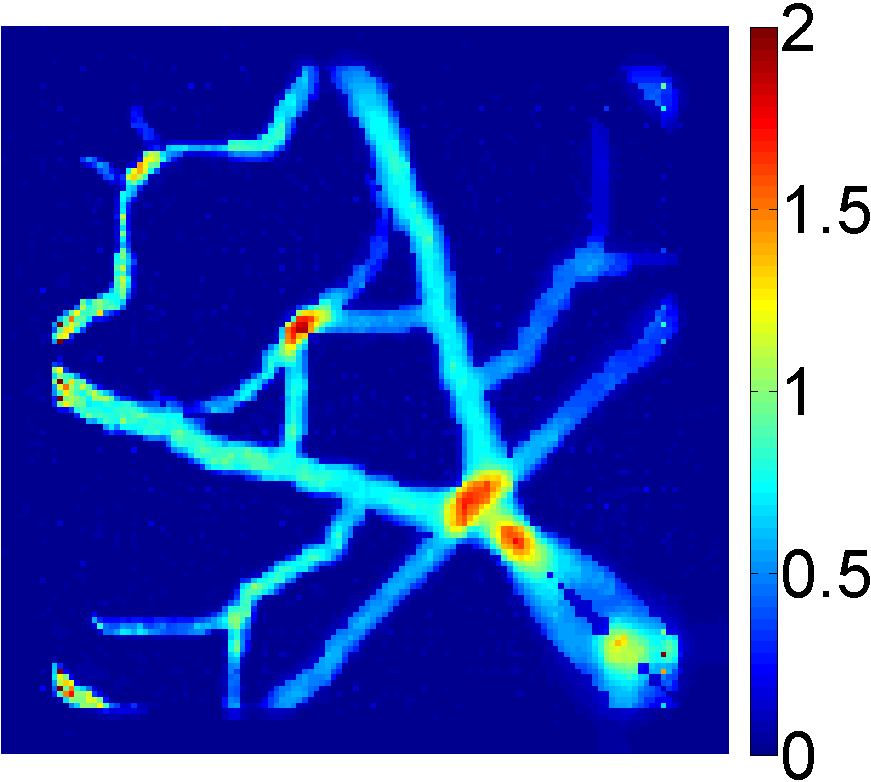}\label{F10b}}
\hspace{0.5cm}
\subfigure[]{\includegraphics[scale=0.13]{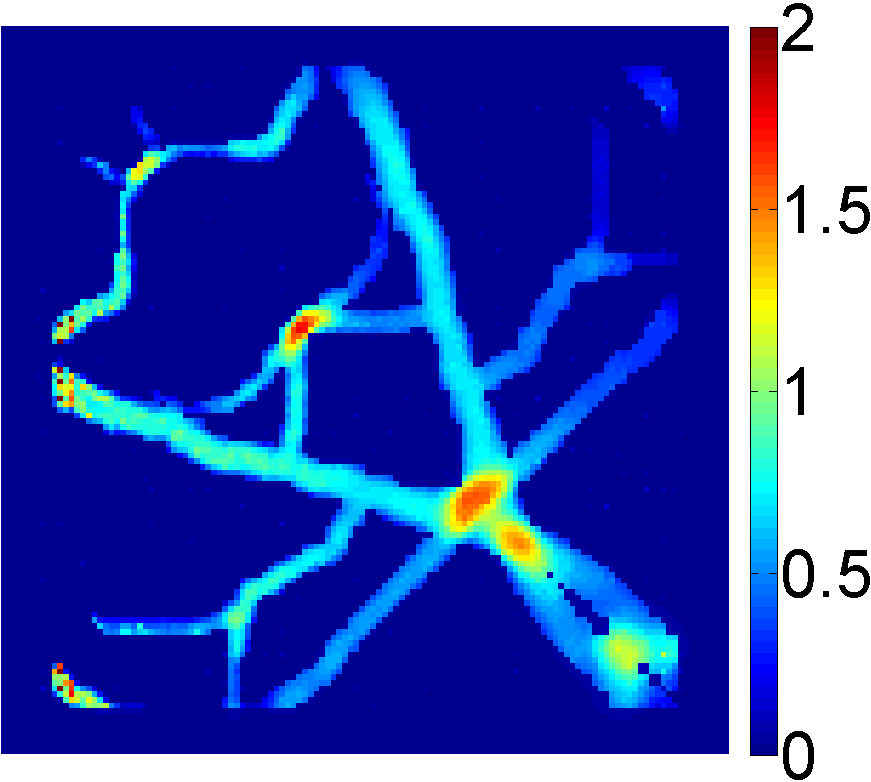}\label{F10c}}
\hspace{0.5cm}
 \subfigure[]{\includegraphics[scale=0.13]{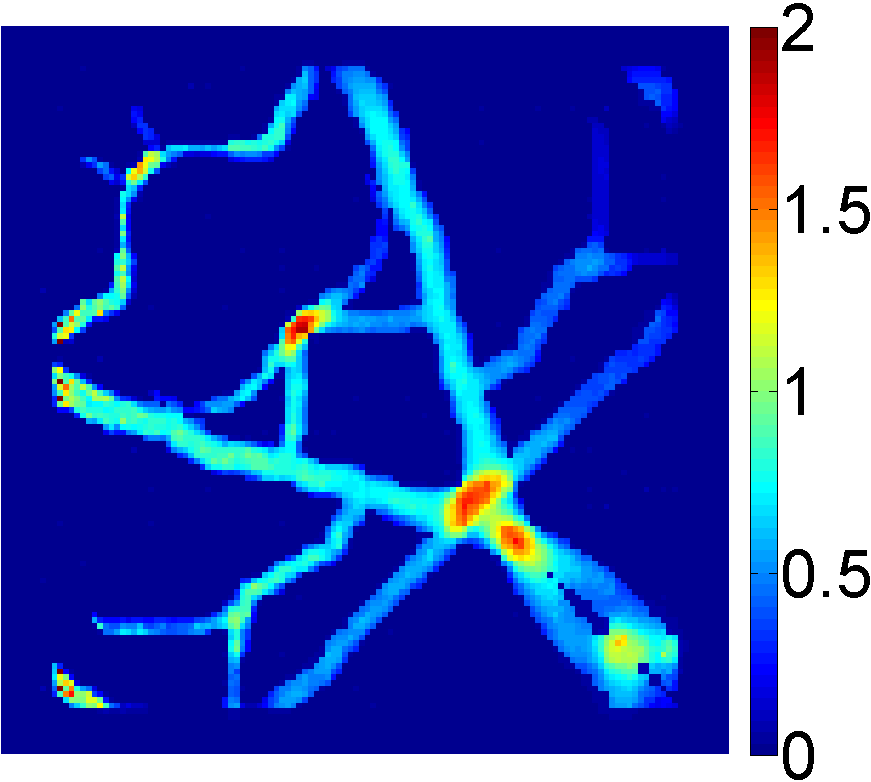}\label{F10d}}
}
\caption{Visualization of 3D images reconstructed by FISTA. (a) Iteration 6 on a fixed grid ($5.34 \times 10^3$s) (b) Iteration 4 on a two-level grid ($3.84 \times 10^3$s) (c) Iteration 16 on a fixed grid ($1.42 \times 10^4$s) (d) Iteration 9 on a two-level grid ($8.59 \times 10^3$s).}
\end{figure}

Figs. \ref{F10a} and \ref{F10c} show images reconstructed by FISTA on the fixed grid at iterations $6$ and $16$ (stopping point), respectively. Figs. \ref{F10b} and \ref{F10d} display images reconstructed by MG FISTA at iterations $4$ and $9$ (stopping point). The MIP visualization provided a different scaling for the reconstructed images, compared to the simulated phantom ($p_\text{exact}$). This different scaling as well as small scale noise do not affect the evaluation of a human observer \cite{Arridge-b}. Thus the reconstructed images were rescaled and thresholded before visualization according to \cite{Arridge-b}. This makes the colorbars equal to that of the simulated phantom, and thus simplifies comparison between the images with respect to the simulated phantom. Accordingly, the visualized image $\bar{x}$ is computed in the form
\begin{align}
\providecommand{\norm}[1]{\left\lVert#1\right\rVert}
\bar{x}=\text{thres}\left(2\frac{x}{\norm{x}_\infty}\hspace{0.05cm},\hspace{0.05cm}0.1\right) \label{forty}
\end{align}
where,
\begin{align}
\text{thres}(v,a) =
\begin{cases}
v, \hspace{0.1cm}\text{if} \hspace{0.1cm} v  \geqslant a\\
0, \hspace{0.1cm}\text{else}
\end{cases}.
\end{align}
Here, $a$ is a thresholding parameter, and the factor $2$ in \eqref{forty} accounts for the maximum amplitude of $p_\text{exact}$.

\section{Conclusion}
We proposed a line search MG optimization approach for PAT. Applied on two grids, our numerical results show that this strategy has improved the speed of ISTA $4$ and $2.35$ times, respectively in 2D and 3D scenarios. A better convergence than MG ISTA was reached by MG FISTA, which was $3$ times faster than fixed-grid FISTA in both 2D and 3D cases.

We derived the adjoint from a first order acoustic system of equations that includes the absorption and dispersion. Our method for deriving the adjoint is in contrast to the method used in \cite{Arridge}, where the adjoint has been derived based on second order acoustic wave equation, and thus does not include absorption and dispersion.

Further studies are needed to extend the proposed MG method to other popular acoustic systems of equations, or other models for describing absorption and dispersion, e.g., \cite{Kowar-b}. The forward implementation of these acoustic systems is often more expensive than the forward model used in our study, and thus solving the corresponding iterative algorithms in a multi-grid setting can be very useful. Additionally, a direct method for quantitative PAT has recently received much attention, where the forward model is treated as a coupled acoustic and optical model, and thus directly links optical properties of medium to time series of acoustic boundary measurements \cite{Haltmeier-b}. Quantitative PAT requires the joint reconstruction of optical absorption and scattering coefficients. This makes the corresponding iterative reconstruction very expensive. It is hoped that an extension of the proposed MG method to the direct quantitative PAT can be useful for improving the reconstruction speed.



%


\bibliography{my_refs}
\bibliographystyle{ieeetran}

\end{document}